\documentclass[aos,seceqn,citesort,dvips]{arximspdf}
\usepackage{graphics}
\doi{10.1214/09-AOS783}
\volume{38}
\issue{4}
\pubyear{2010}
\firstpage{2465}
\lastpage{2498}

\makeatletter

\def\textcyr#1{{\fontencoding{U}\fontfamily{wncyr}\fontsize{11}{11}\selectfont #1}\language4}
\def\textcyrr#1{{\fontencoding{U}\fontfamily{wncyr}\fontsize{8.36}{10.36}\selectfont #1}\language4}

\newcommand{\cha}{\textcyr{ch}}
\newcommand{\chaa}{\textcyrr{ch}}

\newcommand{\dint}{\int\!\!\int}

\newproclaim{Prob}{Problem}[section]
\newtheorem{Th}[Prob]{Theorem}
\newproclaim{Rm}[Prob]{Remark}
\newtheorem{Cor}[Prob]{Corollary}
\newtheorem{Lem}[Prob]{Lemma}

\makeatother

\begin{document}
\begin{frontmatter}

\title{M\"obius deconvolution on the hyperbolic plane with
application to impedance density estimation}
\runtitle{M\"obius deconvolution}

\begin{aug}
\author[A]{\fnms{Stephan F.} \snm{Huckemann}\corref{}\thanksref{t1}\ead[label=e1]{huckeman@math.uni-goettingen.de}},
\author[B]{\fnms{Peter T.} \snm{Kim}\thanksref{t2}\ead[label=e2]{pkim@uoguelph.ca}},\\
\author[C]{\fnms{Ja-Yong} \snm{Koo}\thanksref{t3}\ead[label=e3]{jykoo@korea.ac.kr}} and
\author[A]{\fnms{Axel} \snm{Munk}\thanksref{t4}\ead[label=e4]{munk@math.uni-goettingen.de}}
\runauthor{Huckemann, Kim, Koo and Munk}
\affiliation{Georg-August-Universit\"{a}t G\"{o}ttingen,
University of Guelph, Korea University\break and
Georg-August-Universit\"{a}t G\"{o}ttingen}
\dedicated{Dedicated to Dennis M. Healy Jr. 1957--2009}
\address[A]{S. F. Huckemann\\
A. Munk\\
Institut f\"{u}r Mathematische Stochastik\\
Georg-August-Universit\"{a}t G\"{o}ttingen\\
Goldschmidtstr. 7\\
D-37077 G\"{o}ttingen\\
Germany\\
\printead{e1}\\
\phantom{E-mail: }\printead*{e4}}
\address[B]{P. T. Kim\\
Department of Mathematics\\
\quad and Statistics\\
University of Guelph\\
Guelph, Ontario N1G 2W1\\
Canada\\
\printead{e2}}
\address[C]{J.-Y. Koo\\
Department of Statistics\\
Korea University\\
Anam-Dong Sungbuk-Ku\\
Seoul 136-701\\
Korea\\
\printead{e3}}
\end{aug}

\thankstext{t1}{Supported by DFG Grant MU 1230/10-1 and GRK 1023.}
\thankstext{t2}{Supported by NSERC Discovery Grant 46204.}
\thankstext{t3}{Supported by National Research Foundation of Korea Grant
2009-0075827.}
\thankstext{t4}{Supported by DFG Grant SFB 755 and FOR 916.}

\received{\smonth{10} \syear{2009}}
\revised{\smonth{12} \syear{2009}}

%
\begin{abstract}
In this paper we consider a novel statistical inverse problem on the
Poincar\'{e}, or Lobachevsky, upper (complex) half plane. Here the
Riemannian structure is hyperbolic and a transitive group action comes
from the space of $2\times2$ real matrices of determinant one
via M\"obius transformations. Our approach is based on a deconvolution
technique which relies on the Helgason--Fourier calculus adapted to this
hyperbolic space. This gives a minimax nonparametric density estimator
of a hyperbolic density that is corrupted by a random M\"obius
transform. A motivation for this work comes from the reconstruction of
impedances of capacitors where the above scenario on the Poincar\'{e}
plane exactly describes the physical system that is of statistical
interest.
\end{abstract}

%

\end{frontmatter}

\section{Introduction}\label{intro-scn}

The recovery of objects, for example, densities and functionals
thereof, based on noisy indirect observations, otherwise known as
statistical inverse problems
(see, e.g., \cite{JoSi91}), is scientifically of intense interest. The
literature is vast and we mention only a few selected papers. Most of
the work is concerned with deconvolution on Euclidean spaces. Prominent
approaches are based on wavelet and wavelet-vagulette expansions (see
\cite{AS98,KM03} and \cite{PS09}) or, on singular value
decompositions (see \cite{KPP07,BHMR07} and \cite{MR96}), where
for the latter block thresholding techniques lead to adaptive
estimators (see \cite{CGLT03}). Minimax rates in deconvolution have
been investigated by Fan \cite{Fan91}, and many others. Recently, oracle
inequalities have been proved in \cite{CGPT02}. The specific problem of
boxcar deconvolution and its link to Diophantine approximation have
been investigated in \cite{JKPR04} and~\cite{JR04}. Other methods
include the linear functional strategy (see \cite{Gold99} and variants
thereof). Very popular are Fourier series estimators which have been
known for a long time, and often they become particularly simple
because on Euclidean spaces they can be treated with kernel methods
(see \cite{Fan91,HQ05} and \cite{VEGS08}). We note that even
though our approach also utilizes Fourier methods on groups, due to the
hyperbolic geometry, the resulting estimator cannot be treated by
kernel methods thus complicating our endeavor considerably.

In this paper we provide a novel methodology for statistical object
recovery on the Poincar\'{e} upper half plane which we call the problem
of M\"obius deconvolution. Here the group of M\"obius transformations
is given by all fractions of the form $(az+b)/(cz+d)$ for a complex
number $z$ and a $2 \times2$ matrix with real entries $a,b,c,d$ of
determinant 1. The metric which is invariant under these
transformations is the hyperbolic metric (to be specified in the next
section) which will replace the Euclidean metric in a natural way. In
fact, a key observation is that
this problem can be tackled by generalized Fourier methods similar in
spirit to the Euclidean case. The development of the theory builds on
the foundations laid out in \cite{T85}, Chapter 3. However, extending
Euclidean arguments to this manifold is challenging since the
hyperbolic space is noncompact, the (hyperbolic) geometry is
non-Euclidean and the group of M\"obius transformations is noncommutative.
A~fundamental technical difficulty comes from a lack of a dilation
property that does not extend over
from the Euclidean case.
Despite these difficulties, a generic element is the Riemannian
structure on a manifold,
on which a Laplacian can then be defined. This together with the
generic Euclidean approach as outlined in \cite{MR96} will be the
foundation to what will be presented below.

In addition to the theoretical interest of this novel scenario of M\"
obius deconvolution,
there are important practical applications as well.
One particular situation occurs
in alternating current circuit analysis and design
whenever signals travel through circuit elements as well as whenever
geometries of waveguides change. A simple example of the latter is a
connector to a coaxial cable, say.
Here, the so-called ``reflections'' are modeled on the complex unit disk
and the corresponding ``impedances''
occur in a complex half plane. Usually, these reflections or impedances
are not directly visible but observed through other electrical devices,
such as a ``two-port'' which in turn is modeled by M\"obius
transformations. In particular, a class of so-called ``lossless''
two-ports can be identified with $2 \times2$ matrices of determinant 1.
One particular aspect of M\"obius deconvolution is related to the
temporal decay of impedances of capacitors whereby the above scenario
on the Poincar\'{e} plane exactly
describes the physical system that is of statistical interest.
Other applications include the field of electrical impedance
spectroscopy, as well as electrical impedance tomography. In the
former, measuring varying impedances due to variable ion transport
through biological membranes is currently of high interest in view of
pharmaceutical drug design (see
\cite{EMCG02} as well as
\cite{RS04}). In the latter, in a noninvasive and radiation-free way,
medical imaging can be cost effectively accomplished by measuring
skin-impedances 
(cf. \cite{B01}). Indeed, for successful reconstruction, control of
various errors is of paramount importance (see
\cite{HJDH08}). There is also work in higher-dimensional hyperbolic
spaces with respect to medical imaging (see, e.g., \cite{KR02}).

In statistics there is also some recognition of the Poincar\'{e} plane
and its hyperbolic geometry particularly so because the parameter space
of the Gaussian distribution (with unknown mean and standard deviation)
is this space. Furthermore, it has been shown by several authors that
the Riemannian metric derived from the Fisher information is exactly
hyperbolic (see
\cite{KV97,LMR00} for details).
Obviously, location and dispersion parameters of arbitrary
distributions and random estimators thereof can be viewed within the
Poincar\'{e} plane. Curiously here, the family of Cauchy distributions
play a specific role as being equivariant under M\"obius
transformations (see \cite{M92,M96} for this and its consequences for
parameter estimation).
Based on the above, techniques from the hyperbolic geometry of the
Poincar\'{e} plane are developed
exclusively from a parametric point of view (cf. \cite{MM04}). As far
as the authors are aware, our contribution is the first attempt at
nonparametric developments.


We now summarize the paper. Section \ref{Prelim:scn} is a preliminary
section which introduces
the notation along with the Helgason--Fourier analysis needed for this paper.
Following this, Section \ref{sec:main} presents the main results. In
Section \ref{simulation-scn} we focus on computational aspects of
M\"obius deconvolution illustrating the ideas through simulations. To this
end we introduce in addition to the hyperbolic Gauss the hyperbolic
Laplace distribution.
In Section \ref{engineering-scn} we go into explicit detail with
respect to the M\"obius deconvolution problem for statistically
recovering the temporal decay of impedances of capacitors as outlined
two paragraphs above.
We will briefly sketch the background; however, if
the reader is well versed in this field, then one can start from
Section \ref{example-scn} where we examine a data set that was acquired
through collaboration with the
University of Applied Sciences (Fulda, Germany) that depicts the
physical system of this paper. In particular, we are able to identify
random impedances when only their impedances viewed through random
capacitive two-ports are given.
Following this, technical details of the Poincar\'{e} upper half plane
and the proofs of the main theorems are collected in Appendices \ref{upper-half-plane-scn}
and \ref{sec:proofs}.

As usual for two function $g$ and $f$, write $f\asymp g$ if $f(x)=O
(g(x) )$ and $g(x)=O (f(x) )$ for $x\to\infty$ or $x\to0$,
depending on context.

\section{Preliminaries}\label{Prelim:scn}
In the following let ${\mathbb R}$ and ${\mathbb C}$ denote the real
and complex numbers, respectively. Furthermore, the group of real
$2\times2$ matrices of determinant one is denoted by
%
%
\begin{equation}\label{sl2r} \mathbb{S}\mathbb{L}(2,\mathbb R) := \left\{
g=\pmatrix{
a&b\cr c&d } \dvtx a,b,c,d\in\mathbb R, ad -bc =1 \right\}.
\end{equation}
This defines the group of \textit{M\"obius transformations} $M_g \dvtx
{\mathbb C}\to{\mathbb C}$ by setting for each $g \in\mathbb
{S}\mathbb{L}(2,\mathbb R)$,
%
%
\begin{equation}\label{mtrans}
M_g(z):= \frac{az+b}{cz+d},
\end{equation}
where $M_{g} M_{h} = M_{gh}$ for $g, h \in\mathbb{S}\mathbb
{L}(2,\mathbb R)$.
Let
%
%
\begin{equation}\label{upper-half-plane:def}
\mathbb H := \{z\in\mathbb C\dvtx \operatorname{Im}(z) > 0\}
\end{equation}
be the upper half plane where ``$\operatorname{Re}(z)$''
and ``$\operatorname{Im}(z)$'' denote the real and imaginary parts of a
complex number
$z$, respectively. Then for each $g \in\mathbb{S}\mathbb{L}(2,
\mathbb R)$, the M\"obius
transformation $M_g$ is a bijective selfmap of $\mathbb H$. Moreover,
for arbitrary $z, z' \in\mathbb H$ there exists a (in general not
unique) $g\in\mathbb{S}\mathbb{L}(2, \mathbb R)$ such that $z' = M_g(z')$.

\begin{figure}[b]

\includegraphics{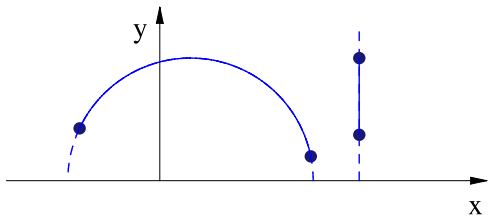}

\caption{The shortest connection (a geodesic segment)
between two points $z,z'\in\mathbb H$ in the hyperbolic geometry is
either a vertical line segment [if $\operatorname{Re}(z)=\operatorname
{Re}(z')$] or an arc on the
circle through $z$ and $z'$ with center on the real axis.}
\label{H_Geods:fig}
\end{figure}

The action of $\mathbb{S}\mathbb{L}(2,\mathbb R)$ on $\mathbb H$,
which is rather
involved, is further discussed in Appendix \ref{upper-half-plane-scn}.
It will be used in the proof of the lower bound in Appendix
\ref{ssec:LowerBound}. For the following we note that M\"obius
transformations preserve the family of vertical lines and circles
centered at the real axis
(cf. Figure \ref{H_Geods:fig}). This is a consequence of the fact, that
M\"obius transformations leave the
cross ratio
\[
c(z_1,z_2,w_1,w_2) = \frac{(z_1-w_1)(z_2-w_2)}{(z_1-z_2)(w_1-w_2)}
\]
invariant. For a detailed introduction (cf. Nevanlinna and Paatero
\cite{NePa64}, Chapter 3).

%
%

The deconvolution, or statistical inverse problem of reconstructing the
density of a random object $X$
on $\mathbb H$, of which we only see a version $Y$ corrupted by an
independent random error $\varepsilon$
on $\mathbb{S}\mathbb{L}(2,\mathbb R)$ can now be formulated as
%
%
\begin{equation}\label{model}
Y = M_\varepsilon(X) .
\end{equation}
A natural geometry for (\ref{model}) is the given by the \textit
{hyperbolic distance} on $\mathbb H$
\[
d(z,z') = \log\frac{1+\sqrt{|c(z,\overline{z'},z',\overline
{z})|}}{1-\sqrt{|c(z,\overline{z'},z',\overline{z})|}}
\]
since for this distance the space of isometries of $\mathbb H$ is
precisely the group of M\"obius transformations, meaning that $d
(M_g(z),M_g(z') )=d(z,z')$ for all $g\in\mathbb{S}\mathbb
{L}(2,\mathbb R)$ and $z,z'
\in\mathbb H$. Here, $\overline{z}= x-iy$ denotes the complex conjugate
of $z=x+iy$. The corresponding \textit{hyperbolic measure} is chosen
such that:
\begin{longlist}
\item its area element $dz$ agrees with the area element $dx\,dy$ of
Lebesgue measure at $z=i$, and
\item it is invariant under M\"obius transformations.
\end{longlist}
In consequence, the Radon--Nikodym derivative of the hyperbolic area
element with respect to the Lebesgue area element at $w=u+iv$ is given
by $v^{-2}$ which is the determinant of the Jacobian
\[
\pmatrix{u_x&u_y\cr v_x&v_y}
\]
of a M\"obius transformation $M$ yielding $M(i) = (ai+b)/(ci+d)
= u+iv$. This can be verified with the complex derivative $M'(i) =
(ci+d)^2 = u_x + iv_x$ and the Cauchy differential equations $u_y=-v_x,
v_y=u_x$; cf. Terras \cite{T85}, Chapter III.

At $z=x+iy \in\mathbb H$ we have hence the hyperbolic area element
%
%
\begin{equation}\label{hyperbolic-measure}
dz := \frac{dx \,dy}{y^2} .
\end{equation}
In addition, in order to properly define below in (\ref{conv_model}) a
convolution of a density on $\mathbb H$ with a density on $\mathbb
{S}\mathbb{L}
(2,\mathbb R)$, a compatible bi-invariant Haar measure $dg$ on $\mathbb
{S}\mathbb{L}
(2,\mathbb R)$ is chosen in Appendix \ref{upper-half-plane-scn}.

Hence, $X$ and $Y$ are random complex numbers in the upper half plane
$\mathbb H$ equipped with the hyperbolic geometry, and $\varepsilon$ is
a random isometric self-map of $\mathbb H$ applied to $X$
by (\ref{mtrans}). The problem of the \textit{M\"obius deconvolution}
can be made precise as follows. A density on the upper complex half
plane with respect to
the hyperbolic measure is called a hyperbolic density. Densities on
$\mathbb{S}\mathbb{L}(2,\mathbb R)$ are taken with respect to the
Haar measure $dg$.
\begin{Prob}\label{prob_statement} Under the model (\ref{model}) estimate
nonparametrically, the hyperbolic density $f_X$ of $X$
from the hyperbolic density $f_Y$ of $Y$ when the density
$f_{\varepsilon}$ on $\mathbb{S}\mathbb{L}(2,\mathbb R)$ is known.
\end{Prob}

We note that this setup assumes underlying i.i.d. $X_1,\ldots,X_n$
corrupted by i.i.d. errors $\varepsilon_1,\ldots,\varepsilon_n$, also
independent of $X_j$ $(j=1,\ldots,n)$ giving observations $Y_j =
M_{\epsilon_j}(X_j),j=1,\ldots,n$.\vspace*{1pt}

We will base our work using Fourier or singular value decomposition
methods that are common
for the Euclidean case (see \cite{KPP07,BHMR07} and
\cite{MR96}) and the fact
that the densities of (\ref{model}) are related by the convolution
%
%
\begin{equation}\label{conv_model}
f_Y(z) = (f_\varepsilon*f_X)(z) =: \int_{\mathbb{S}\mathbb
{L}(2,\mathbb R)}
f_{\varepsilon}(g)
f_X (M_{g^{-1}}(z) ) \,dg.
\end{equation}
This fact is a consequence of (\ref{conv_lift:def}) in Appendix
\ref{upper-half-plane-scn}.

From here on, we will make the abbreviation $\mathbb{S}\mathbb
{L}(2):= {\mathbb{S}\mathbb{L}(2,\mathbb
R)}$, as well as to write (\ref{mtrans}) as simply $g(z)$ or $gz$ for
$g \in\mathbb{S}\mathbb{L}(2)$ and $z \in\mathbb H$ whenever the
context is clear.

\subsection{\texorpdfstring{Fourier analysis on the Poincar\'{e}
plane}{Fourier analysis on the Poincare plane}}
\label{fourier-analysis-scn}
For purposes of M\"obius deconvolution for Problem \ref
{prob_statement} we sketch the {Helgason formulation} of hyperbolic
Fourier calculus which can be found in more detail in \cite{T85},
Chapter 3.2.
The \textit{Helgason--Fourier} transform of $f \in C_c^{\infty
}(\mathbb
H)$ with the latter being the space of real valued functions with
compact support in $\mathbb H$ with derivatives of all orders, is
defined as the function
\[
\mathcal{H}f(s,k) := \int_{\mathbb H} f(z)
(\operatorname{Im}(k(z) ) )^{\overline{s}} \,dz
\]
analytic for $(s,k) \in\mathbb C \times\mathbb{S}\mathbb{O}(2)$
where overline
denotes complex
conjugation. Here,
\[
k = k_u = \pmatrix{\cos u&\sin u\cr-\sin u&\cos u}
\in\mathbb{S}\mathbb{O}(2) \subset\mathbb{S}\mathbb{L}(2)
\]
is naturally identified with $u\in[0,2\pi)$ acting on $\mathbb H$ as
the M\"obius transformation
\[
M_k(z) =\frac{z\cos u + \sin u}{\cos u - z\sin u}
\]
as defined in (\ref{mtrans}) (cf. Appendix \ref{upper-half-plane-scn}).
Note that for all $s \in\mathbb C$, $z =x+iy\to y^s$ and $z\to
(\operatorname{Im}(k(z) ) )^s$ are eigenfunctions with corresponding
eigenvalues $s(s-1)$ of the Laplace--Beltrami operator
%
%
\begin{equation}
\label{laplace}
\Delta= y^2 \biggl(\frac{\partial^2}{\partial x^2} +
\frac{\partial^2}{\partial y^2} \biggr)
\end{equation}
on $\mathbb H$. With the \textit{spectral measure}
\[
d\tau= \frac{1}{8\pi^2} t \tanh(\pi t) \,dt \,du
\]
on $\mathbb R \times
\mathbb{S}\mathbb{O}(2)$ the \textit{inverse Helgason--Fourier}
transform is given by
%
%
\begin{equation}\label{inv_transfo}
f(z)
=\int_{t\in\mathbb
R}\int_{u=0}^{u=2\pi} \mathcal{ H}f \biggl(\frac{1}{2} +it,k_u \biggr)
(\operatorname{Im}(k(z) ) )^{{1}/{2}+it} \,d\tau,
\end{equation}
where $i^2 =-1$.
The following result justifies these definitions: mapping to the
Helgason--Fourier transform extends to an isometry $L^2(\mathbb H,dz)
\to L^2(\mathbb R\times\mathbb{S}\mathbb{O}(2),d\tau)$; that is, we
have the \textit{Plancherel identity}
%
%
\begin{equation}\label{planch_eq}
\int_{\mathbb H} |f(z)|^2 \,dz = \int_{t\in\mathbb
R}\int_{u=0}^{u=2\pi} \biggl|\mathcal{ H}f \biggl({\frac1 2}+it,k_u
\biggr) \biggr|^2 \,d\tau,
\end{equation}
where we denote the space of square integrable functions over some space
by~$L^2$. We note that
$f\in L^2(\mathbb H,dz)$ is $\mathbb{S}\mathbb{O}(2)$-invariant if
and only if
$\mathcal{H}f\in L^2(\mathbb R\times\mathbb{S}\mathbb{O}(2),d\tau
)$ is $\mathbb{S}\mathbb{O}(2)$-invariant.
Thus, for numerical computations (\cite{T85}, pages 141 and 149)
for an $\mathbb{S}\mathbb{O}(2)$-invariant function $f$, transforms
and inverse
transforms can be considerably simplified
%
%
\begin{eqnarray}\label{transfo_invariant}\quad
\mathcal{H}f \biggl(\frac{1}{2} +it \biggr) &=& 2\pi\int
_0^{\infty} f(e^{-r}i) P_{-{1}/{2} +it}(\cosh r) \sinh r \,dr
,\nonumber\\
f(e^{-r}i) &=&\frac{1}{4\pi} \int_{-\infty}^{\infty} \mathcal{H}f
\biggl(\frac{1}{2} +it \biggr) P_{-{1}/{2} +it}(\cosh r) t \tanh
(\pi t) \,dt
\end{eqnarray}
with the {Legendre function}
\[
P_a(c) := \frac{1}{2\pi}\int_0^{2\pi} \bigl(c + \sqrt{c^2-1} \cos
(\phi) \bigr)^{a} \,d\phi.
\]

Throughout this work we will use the following assumptions:
\begin{enumerate}[(D.1)]
\item[(D.1)] all densities are square-integrable
\[
f_X, f_Y \in
L^2(\mathbb H,dz),\qquad f_{\varepsilon} \in L^2(\mathbb{S}\mathbb
{L}(2),dg) ;
\]
\item[(D.2)] the error density $f_{\varepsilon}$ is {bi-invariant}
\[
f_{\varepsilon}(agb) = f(g) \qquad\forall g \in\mathbb
{S}\mathbb{L}(2), a,b\in
\mathbb{S}\mathbb{O}
(2) ;
\]
\item[(D.3)] $f_X \in\mathcal{F}_{\alpha}(Q)$ for a Sobolev
ball
\[
\mathcal{F}_{\alpha}(Q) = \{f\in L^2(\mathbb H,dz)\dvtx \|\Delta^{
{\alpha}/{2}}f\|^2
\leq Q\}
\]
with $\alpha> 1$ and $Q>0$.
\end{enumerate}
Here, $\Delta^{{\alpha}/{2}}f$ denotes the unique function $h\in
L^2(\mathbb H,dz)$ with $\mathcal{H} h(s,k) =
\overline{s(s-1)}{}^{{\alpha}/{2}}\mathcal{H} f(s,k)$.

As detailed in Appendix \ref{upper-half-plane-scn}, the isometry
$\mathbb{S}\mathbb{L}(2)/\mathbb{S}\mathbb{O}(2) \to\mathbb H\dvtx
g\mathbb{S}\mathbb{O}(2) \mapsto
M_g(i)$ preserves the action of $\mathbb{S}\mathbb{L}(2)$. Hence a density
$f_{\varepsilon}$ satisfying (D.2) can
be regarded
as an $\mathbb{S}\mathbb{O}(2)$-invariant
mapping $\mathbb H \to\mathbb R$. In particular, in case of (D.1) and
(D.2), the Helgason--Fourier transform $ \mathcal{ H}f_{\varepsilon}(z)$
is well defined, we have
%
%
\begin{equation}\label{helg_trans_conv}
\mathcal{H}f_Y(s,k) = \mathcal{H}(f_{\varepsilon
}*f_X)(s,k)=\mathcal{ H}f_{\varepsilon}(s)\cdot\mathcal{H}f_X(s,k);
\end{equation}
\cite{T85}, page 149. One final
assumption to be made is the following.
\begin{enumerate}[(D.4)]
\item[(D.4)] $\exists$ constants $\beta,\gamma, C_1,C_2 >0\mbox{
such that }$
\[
C_1 \exp\biggl\{-\frac{|s|^{\beta}}{\gamma} \biggr\} \leq|\mathcal
{H}f_{\varepsilon}(s) |\leq
C_2 \exp\biggl\{-\frac{|s|^{\beta}}{\gamma} \biggr\}\quad \ \  \forall s = \frac
{1}{2}+it,\quad \ \  t\in\mathbb R .
\]
\end{enumerate}
As an example, the hyperbolic Gaussian-distribution (see Section \ref
{sec:main} below) satisfies (D.4).

Of course, any density on the upper half plane (or on the unit disk)
can be rescaled
with respect to hyperbolic measure. One example, using the normalized
squared absolute cross ratio $c(z,\overline{\theta},\overline
{z},\theta
)$, has been kindly provided by one of the referees,
\[
\frac{1}{\pi} \biggl(\frac{|z-\overline{z}||\theta-\overline{\theta
}|}{|z-\overline{\theta}|^2} \biggr)^2\, dz=\frac{4\sigma^2\,dx\, dy}{\pi
((x-\mu)^2+(y+\sigma)^2 )^2},\qquad z=x+iy \in\mathbb H,
\]
with a hyperbolic parameter $\theta=\mu+i\sigma\in\mathbb H$; for
example, one could take $x$ and $y$ as suitable estimators of the
location and dispersion parameters of another distribution (cf. \cite
{M92,M96}). Note that this density is not $\mathbb{S}\mathbb
{L}(2)$-invariant; rather it
is equivariant with respect to the $\mathbb{S}\mathbb{L}(2)$ action
on both variable $z$
and parameter $\theta$. In the context of this research one is
interested also in $\mathbb{S}\mathbb{L}(2)$-invariant densities.
Such can be generated
from suitable densities on $y \in[1,\infty)$. Moreover, $\mathbb
{S}\mathbb{L}
(2)$-invariant functions additionally fulfilling (D.3) can be obtained
by applying the inverse Helgason transform (\ref{transfo_invariant})
to suitable functions on $s = \frac{1}{2}+it$,
$-\infty< t<\infty$.
If the function is even in $t$, then the inverse Helgason transform
thus obtained is real. It is, however, not necessarily nonnegative. As
a consequence of $P_{-1/2+it}(\cosh r) >0$ for $t=0$, nonnegativity
can be obtained if the function tends sufficiently fast to zero as
$t\to
\infty$. Numerical experiments indicate that one may consider for
$\tau
> -1/4$ and $\alpha>1$ a suitable multiple of a power of a Cauchy
density in the spectral domain
%
%
\begin{equation}\label{invH_laplace_dens:eq} \mathcal{H}h_{\alpha
,\tau
}(s,k) \propto \frac{1}{ (\tau-s(s-1) )^{\alpha}} = \frac
{1}{ (\tau+1/4+t^2 )^{\alpha}}
\end{equation}
giving an invariant \textit{hyperbolic Laplace density} $h_{\alpha
,\tau}
\in\mathcal{F}_{\alpha}(Q)$.
This density can be lifted as in Appendix \ref{upper-half-plane-scn}
giving a bi-invariant density $\widetilde{h}_{\alpha,\tau}$ on
$\mathbb{S}\mathbb{L}
(2)$. In particular $\widetilde{h}_{\alpha_1,\tau_1} *h_{\alpha
_2,\tau
_2}\in\mathcal{F}_{\alpha_1+\alpha_2}(Q)$ for $\tau_1,\tau_2>-1/4$,
$\alpha_1,\alpha_2>1$ and suitable $Q>0$.

We will not elaborate further on this topic, but we mention
that motivated by our research and by many potential applications,
the task of generalizing non-Gaussian distributions to hyperbolic
spaces may lead to a new field of challenging research.

\section{Main results}\label{sec:main}

Let us begin with the definition of the Helgason--Fourier transform of
the generalized derivative of the empirical distribution $f^{(n)}_Y(z)
= \frac{1}{n} \sum_{j=1}^n\delta_{Y_j}(z)$ where $Y_1,\ldots, Y_n$
is a random sample in $\mathbb H$
%
%
\begin{equation}\label{emp_transfo}
\mathcal{H}f^{(n)}_Y(s,k) = \frac{1}{n} \sum_{j=1}^n (\operatorname{Im}
(k(Y_j) ) )^{\overline{s}} .
\end{equation}
Obviously
%
%
\begin{equation}\label{emp_transfo_exp_value}
\mathbb{E}\mathcal{H}f^{(n)}_Y(s,k) = \mathcal{H}f_Y(s,k) ,
\end{equation}
where ``$\mathbb{E}$'' denotes expectation.
We estimate the Helgason transform of an $\mathbb{S}\mathbb
{O}(2)$-invariant density as
well by an $\mathbb{S}\mathbb{O}(2)$-invariant estimator
%
%
\begin{equation}\label{emp_transfo_invariant}
\mathcal{H}f^{(n)}_Y(s) :=\mathcal{H}f^{(n)}_Y(s,{\mathbf{I}} )=\frac
{1}{n} \sum_{j=1}^n \operatorname{Im}(Y_j)^{\overline{s}}
\end{equation}
with the identity element $\mathbf{I}\in\mathbb{S}\mathbb{O}(2)$.

From the Helgason--Fourier transform (\ref{emp_transfo})
we build an estimator by using (\ref{helg_trans_conv}) and the inverse
Helgason--Fourier transformation (\ref{inv_transfo}) with a suitable
cutoff $T>0$
%
%
\begin{equation}\label{inv_emp_transfo}
f^{(n,T)}_{X}(z)
:= \int_{|t|<T}\int_{u=0}^{u=2\pi} \frac{\mathcal{ H}f^{n}_Y
({1}/{2} +it,k_u )}{\mathcal{ H}f_{\varepsilon} (
{1}/{2} +it )} (\operatorname{Im}(k_u(z) ) )^{{1}/{2}+it}\,
d\tau
\end{equation}
for the density $f_X$. This is well defined if $\mathcal{
H}f_{\varepsilon} \neq0$ is bounded from below on compact sets which
is guaranteed
under assumption (D.4). Even though we consider in this section the
general case, we note in view of (\ref{emp_transfo_invariant}) and
(\ref{transfo_invariant}) that for the estimation of an $\mathbb
{S}\mathbb{O}
(2)$-invariant density $f_X$, we can use the simpler
\[
f^{(n,T)}_{X}(e^{-r}i)
:= \frac{1}{4\pi} \int_{-T}^T\frac{\mathcal{ H}f^{n}_Y (
{1}/{2} +it )}{\mathcal{ H}f_{\varepsilon} ({1}/{2}
+it )} P_{-{1}/{2}+it}(\cosh r) t \tanh(\pi t) \,dt .
\]

As the first main result we have:
\begin{Th}\label{main-thm}For $f_X, f_Y$ and $f_{\varepsilon}$
satisfying \textup{(D.1)--(D.3)}, and $\mathcal{ H}f_{\varepsilon} \neq
0$ bounded from below on compact sets, there is a
constant $C>0$ not depending on $T,\alpha, Q$ and $n$ such that
\[
\mathbb{E}\bigl\|f_X^{(n,T)} - f_X\bigr\|^2 \leq C \sup_{|t|\leq T} \biggl|\mathcal{
H}f_{\varepsilon} \biggl(\frac{1}{2} +it \biggr) \biggr|^{-2} \frac{T^{2}}{n}
+ QT^{-2\alpha}
\]
as $n \to\infty$.
\end{Th}

If the corruption by error is smooth enough, or equivalently if the
asymptotic rate of
the decay of its Helgason--Fourier transform is suitable, the cutoff
$T$ can be adjusted appropriately
to obtain the following rates.
\begin{Th}\label{opt-T-cor} Suppose that $f_X, f_Y$ and
$f_{\varepsilon
}$ satisfy \textup{(D.1)--(D.4)}.
Then by letting
$T= (\frac{\gamma}{2}\log n - \frac{\eta\gamma}{2} \log(\log
n) )^{{1}/{\beta}}$
where
$\eta\geq2(\alpha+1)/\beta$
\[
\mathbb{E}\bigl\|f_X^{(n,T)} - f_X\bigr\|^2 \leq Q \biggl(\frac{\gamma}2 \log n
\biggr)^{-{2\alpha}/{\beta}}\bigl(1+o(1)\bigr)
\]
as $n \to\infty$ where $\alpha$ is from condition \textup{(D.3)}.
\end{Th}

The optimal rate of a power of $\log n$ in case of error smoothness
(D.4) is in agreement with Euclidean results. On the real line,
condition (D.4) corresponds to {supersmooth} errors for which Fan
\cite{Fan91} establishes the same type of rate. This rate has also been
established by Butucea and Tsybakov \cite{BT07} in case of additionally
supersmooth signals. For a scenario corresponding to our setup on
compact Lie groups, see \cite{KR01} and more general on any compact
manifold, see \cite{KK05}, where similar rates have been found.

The above results are minimax in the sense that the rate of
convergence is matched by a corresponding
lower bound. We have the following theorem.
\begin{Th}\label{thm:LowerBound} Suppose that $f_X, f_Y$ and
$f_{\varepsilon}$ satisfy \textup{(D.1)--(D.4)}. Then for some constant
$C > 0$, we have
\[
\inf\sup\mathbb{E}\|f^n - f_X\|^2 \geq C (\log n)^{-{2\alpha
}/{\beta
}}
\]
as $n \to\infty$, where the infimum is taken over all estimators
$f^n$ and the supremum over all $f_X \in\mathcal{F}_{\alpha}(Q)$.
\end{Th}

Recall the Gaussian density $g_{\rho}$ on the real line with zero mean
and variance $2\rho>0$ can be characterized as yielding the solution of
the heat equation
\[
(\Delta-\partial_{\rho})u =0
\]
with initial condition $u(z,0) = f(z)$ by
%
%
\begin{equation}\label{heat-eqn}
u(z,\rho) = (g_{ \rho}*f)(z) .
\end{equation}
Similarly on $\mathbb H$, the density $g_{\rho}$ giving the solution
of the heat equation by (\ref{heat-eqn}) is also called the Gaussian
density for $\mathbb H$.
Here, the Laplace--Beltrami operator $\Delta$ would be defined by
(\ref
{laplace}), and the
convolution in (\ref{heat-eqn}) would be defined
as in (\ref{conv_model}).
Using (\ref{helg_trans_conv}) for $\mathbb{S}\mathbb
{O}(2)$-invariant $f$ and $u$ it is
easily seen that
$\mathcal{H}g_{ \rho}(s) \propto e^{\overline{s(s-1)} \rho}$.
Consequently, in terms of assumption (D.4), the Gaussian density satisfies
$\beta= 2$ and $\gamma= 1/\rho$. We have the following result.
\begin{Cor}\label{opt-T-cor-1} For $f_X$ and $f_Y$ satisfying
\textup{(D.1)--(D.3)}
consider corruption according to a Gaussian distribution
$f_{\varepsilon
} = g_{\rho}$.
Then by letting
$T^2= \frac{1}{4\rho} [\log n - \eta\log(\log n) ]$
where
$\eta\geq1+ \alpha$,
\[
\mathbb{E}\bigl\|f_X^{(n,T)} - f_X\bigr\|^2 \asymp(\log n )^{-\alpha}
\]
as $n \to\infty$ gives the optimal rate of convergence.
\end{Cor}

\section{Computations and simulations}
\label{simulation-scn}

In this section we elaborate on computational aspects, simulations and,
in particular, discuss the
Gaussian distribution on~${\mathbb H}$. We begin by first discussing
methods for choosing the truncation parameter.

\subsection{Estimating truncation parameter}
\label{cross-valid-scn}
A popular technique for data-driven choice of a truncation parameter is
least squares cross-validation (see 
\cite{DMR96} or
\cite{W94}, Chapter 3.3). We will discuss how that technique can be
adapted to our setting.
For a given random sample $Y_1,\ldots, Y_n$, an optimal cutoff
$T=T^*_n>0$ minimizes the mean integrated squared error
\[
T^*_n = \mathop{\arg\min}_{T>0} \biggl\{\mathbb{E}\biggl(\int_{\mathbb H}
\bigl(f_{X}^{(n,T)}(z) \bigr)^2 \,dz \biggr)-2\mathbb{E}\biggl(\int_{\mathbb H} f_X(z)
f_{X}^{(n,T)}(z) \,dz \biggr) \biggr\} .
\]
Instead of deriving a minimizer of the above we content ourselves with
minimizing a suitable estimator.
Obviously, $\int_{\mathbb H} (f_{X}^{(n,T)}(z) )^2 \,dz$ is an
unbiased estimator of the
first term. Let
\[
f_{X}^{(n,T,l)} := \int_{-T}^T \int_0^{2\pi}\frac{1}{n-1} \sum
_{j\neq
l}\operatorname{Im}(k(Y_j) )^{{1}/{2}-it} \operatorname
{Im}(k(Y_l) )^{
{1}/{2}+it}\frac{d\tau}{\mathcal{H}f_{\varepsilon} (
{1}/{2}+it )}
\]
and therefore choose
\[
T_n := \mathop{\arg\min}_{T>0} \Biggl(\int_{\mathbb H}
\bigl(f_{X}^{(n,T)}(z) \bigr)^2 \,dz - \frac{2}{n} \sum_{l=1}^n
f_{X}^{(n,T,l)} \Biggr),
\]
which is an estimate for an optimal $T=T^*_n$.

Alternatively, we can use the result of Corollary \ref{opt-T-cor} and set
\[
T= \biggl[\frac{\gamma}{2}\log n - \frac{\gamma}{2} \log(\log n)^\eta
\biggr]^{{1}/{\beta}} .
\]

We are aware of the fact that cross-validation in general suffers from
too large variability and, of course, more involved parameter selection
methods could be generalized here as well (see, e.g.,
\cite{DG04,DHM08,N08,PS05} and \cite{S03} among many others).
However, we do not pursue this issue any further in this paper.

\subsection{Simulation of the Gaussian distribution}
For simulation we use the analog $g_\rho$ of the Gaussian distribution
on the upper half plane introduced above. Recall that
by a more subtle argument (see \cite{T85}, pages 153 and 155), the
inverse transform is obtained in polar coordinates [for any $k\in
\mathbb{S}\mathbb{O}(2)$]
\begin{eqnarray*}
g_{ \rho} (k(e^{-r}i) )&=&g_{ \rho}(e^{-r}i)
= \frac{1}{\sqrt{4\pi\rho}^3}\sqrt{2}e^{- \rho/4}\int
_r^{\infty
}\frac{b e^{-b^2/4 \rho}\,db}{\sqrt{\cosh b -\cosh r}}\\
&=&\!:\frac{\widetilde{g}_{ \rho}(r)}{2\pi\sinh r} ;
\end{eqnarray*}
that is, for a $\rho_X$-Gaussian distributed $\mathbb{S}\mathbb
{O}(2)$-invariant random
object $X$ on $\mathbb H$ and an $\mathbb{S}\mathbb{O}(2)$-invariant
subset $A \subset
\mathbb H$,
\[
{\mathbb P}\{X \in A\} = \int_{\{r\geq0\dvtx e^{-r}i\in A\cap\mathbb H\}}
\widetilde{g}_{ \rho_X}(r)\,
dr
\]
[see (\ref{hyperbolic-measure-polar-coordinates})]. Hence, in order to
simulate $X$ from an invariant Gaussian distribution we simulate $r_X
\sim\widetilde{g}_{\rho_X}$ on $\mathbb R$ and $u_X$ uniform on
$[0,2\pi)$; then
\[
X = k_{u_X} \circ R_{r_X} (i) = k_{u_X}(e^{-r_X}i) .
\]
Note that $g_{\rho_X} \in\mathcal{F}_{\alpha}(Q)$ for all $\alpha>0$
with suitable $Q=Q_{\alpha}>0$. That is, $f_X$ satisfies condition (D.3).

Similarly, in order to simulate $\varepsilon$ from a bi-invariant
$\rho
_{\varepsilon}$-Gaussian distribution on $\mathbb{S}\mathbb{L}(2)$
we consider
$r_{\varepsilon}\sim\widetilde{g}_{\rho_{\varepsilon}}$ and
$k_{u_{\varepsilon}},k_{u'_{\varepsilon}}$ independent and uniform on
$[0,2\pi)$. Then
\[
Y = k_{u_{\varepsilon}}\circ R_{r_{\varepsilon}} \circ
k_{u'_{\varepsilon}} (X)
\]
[see (\ref{model}) and (\ref{SL(2)-polar-decomposition})].

According to Theorem \ref{main-thm}, (\ref{emp_transfo}) and (\ref
{inv_emp_transfo}) we can then estimate the density $f_X$ by
\[
f^{(n,T)}_{X}(z)
= \int_{|t|<T}\int_{u=0}^{u=2\pi} \frac{{1}/{n} \sum_{j=1}^n
(\operatorname{Im}(k_u(Y_j) ) )^{{1}/{2} -it}}
{e^{- (t^2+{1}/{4} ) \rho_{\varepsilon}}} (\operatorname{Im}
(k_u(z) ) )^{{1}/{2}+it} \,d\tau.
\]

By $\mathbb{S}\mathbb{O}(2)$-invariance it suffices to estimate for
$z = e^{-r}i$
only, hence we estimate $\widetilde{g}_{ \rho_X}(r)$ by
\[
\widetilde{f}^{(n,T)}_X(r) := 2\pi\sinh r f^{(n,T)}_{X}(e^{-r}i)
\]
with the integral simplified as in (\ref{transfo_invariant}).

%
%
\begin{figure}

\includegraphics{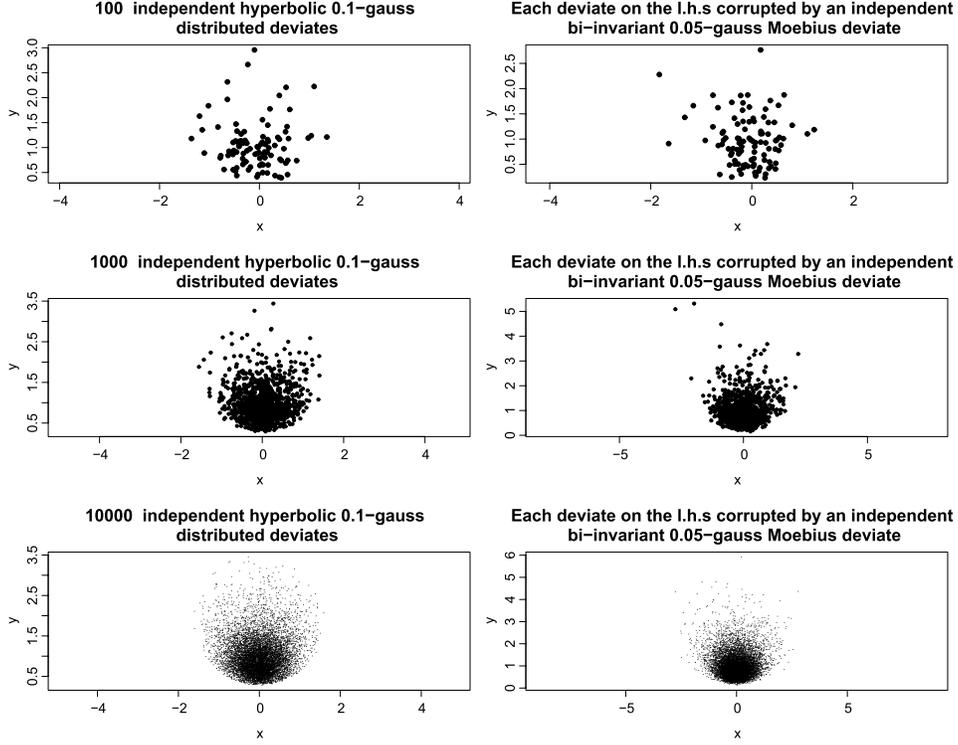}

\caption{Three simulated data samples of $n=100$
(top row), $n=1000$ (middle row) and $n=10\mbox{,}000$ (bottom row)
independent random data points on the upper half plane. Left
row: original independent invariant $\rho_X$-Gaussian
distributed data points, right row: transformed data points
under independent n bi-invariant $\rho_\varepsilon$-Gaussian
distributed $\mathbb{S}\mathbb{L}(2)$ transformations
(right).}\label{sim-pts-fig}
\end{figure}

%
%
\begin{figure}

\includegraphics{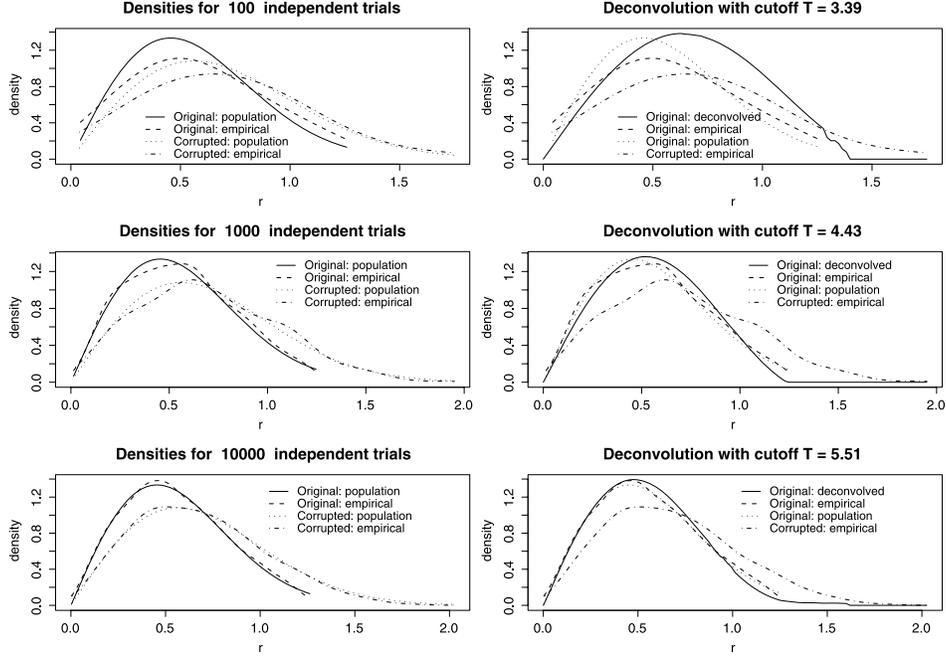}

\caption{Left-hand side: population and empirical
densities times hyperbolic area [corresponding to
$\widetilde{g}_{\rho_X}(r)$ and
$\widetilde{g}_{\rho_X+\rho_{\varepsilon}}(r)$] along the first
polar coordinate $r$ of the data depicted and described in Figure
\protect\ref{sim-pts-fig}. Top row: $n=100$, middle row: $n=1000$ and bottom
row: $n=10\mbox{,}000$. Right-hand side: additionally the respective estimate
times hyperbolic area: $\widetilde{g}{}^{(n,T)}_X(r)$ of the original
density by M\"{o}bius deconvolution. The corresponding optimal cutoff
parameters $T=T_n$ have been estimated by least squares
cross-validation as in Section
\protect\ref{cross-valid-scn}.}
\label{deconv_fig}
\end{figure}

In the following simulation we consider an original
distribution with $\rho_X = 0.1$ under a corrupting M\"obius transformation
distributed with $\rho_{\varepsilon}=0.05$. From this we create three
data sets with different sample sizes: $n=100,1000$ and 10,000.
Figure \ref{sim-pts-fig}
shows the original $X_1,\ldots, X_n$ and the corrupted data $Y_1 =
M_{\varepsilon_1} (X_1),\ldots,Y_n=M_{\varepsilon_n} (X_n)$ in
cartesian coordinates in the upper half plane.
In Figure \ref{deconv_fig}
we show the corresponding densities times hyperbolic area on $[0,\infty
)$. Note that these are then densities in the usual sense; that is,
their integrals with respect to Lebesgue measure on $[0,\infty)$ are 1.
The density estimation by deconvolution has been obtained from the
observed data $Y_1,\ldots, Y_n$ by the proposed method. For the
deconvolution, since only the optimal rate
\[
T \approx\biggl(\frac{1}{4\rho_{\varepsilon}} \log n \biggr)^{
{1}/{4}} =\cases{2.19, &\quad $(n=100)$,\cr
2.42, &\quad $(n=1000)$,\cr
2.61, &\quad $(n=10\mbox{,}000)$,
}
\]
is guaranteed by Corollary \ref{opt-T-cor-1}, we have used the estimate
{via} least squares cross-validation as proposed in Section \ref
{cross-valid-scn}.

\subsection{Simulation of the hyperbolic Laplace distribution}

Using formula (\ref{transfo_invariant}) directly with (\ref
{invH_laplace_dens:eq}) to simulate Laplace $(\alpha,\tau)$-deviates
for Laplace distributed $\mathbb{S}\mathbb{L}(2)$ error corruption,
we obtain results
similar to the ones reported above. Due to the oscillation of the
Legendre polynomials to be evaluated, however, the computational time
is much longer. In analogy to Theorem \ref{opt-T-cor} we have the upper
bound $O(n^{- {\alpha}/({1+\alpha})})$ for the choice $T =
n^{{1}/({2(\alpha+1)})}$.

\section{Impedance density estimation in AC driven circuits}
\label{engineering-scn}

For the convenience of the reader, we begin this section with a review
of classical electrical engineering theory specifically tailored to the
application of hyperbolic statistics in mind. For the underlying
engineering terminology we refer to standard textbooks such as 
\cite{Collin92}. More mathematical approaches are explained in 
\cite{H82,AH03} and~\cite{T85}, Chapter 3. 
In the following we rephrase this problem in the language
of statistics. We are then able to identify a typical problem as a
novel inverse problem in hyperbolic space.

Here, general M\"obius transformations appear with complex
coefficients $a$, $b$, $c$, $d$ in (\ref{sl2r}). Moreover, hyperbolic space
materializes in the form of the upper half plane $\mathbb H$, the open
{unit disk}
$\mathbb D:=\{w\in\mathbb C\dvtx |w| <1\}$, and the open right half plane
$-i\mathbb H := \{\zeta\in\mathbb C\dvtx \operatorname{Re}(\zeta) > 0\}
$. With the
notation of (\ref{upper-half-plane:def}), all are related to one
another by M\"obius transformations, the first is usually called the
\textit{Cayley transform}
%
%
\begin{eqnarray}
\label{Cayley-transform}
w &=& \mathcal{C}(z):=\frac{z-i}{z+i},\qquad z=i \frac{1+w}{1-w},\qquad
i\zeta=z,\nonumber\\[-8pt]\\[-8pt]
w&=&\frac{\zeta-1}{\zeta+1},\qquad\zeta= \frac{1+w}{1-w} .\nonumber
\end{eqnarray}

\subsection{Complex impedance in AC circuits}\label{Complex_imp_scn}
\begin{figure}
\begin{tabular}{cc}

\includegraphics{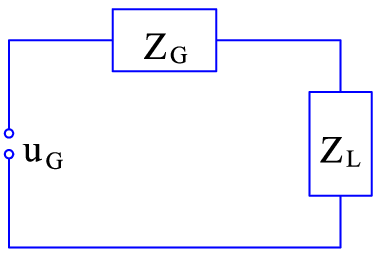}
 & \includegraphics{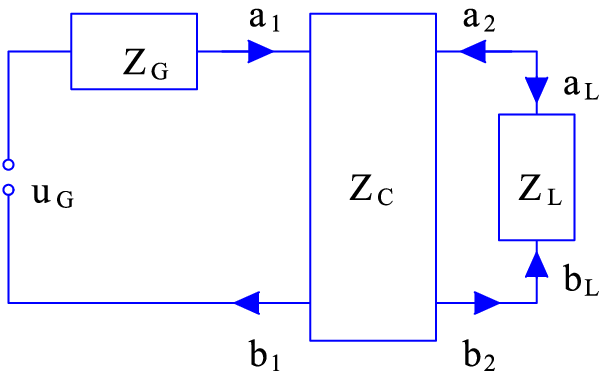}\\
(a) & (b)
\end{tabular}
\caption{Basic circuit models
for signal processing. Left: one-port, right: two-port. \textup{(a)} Serial circuit
with generated voltage $u_G$, generator impedance $Z_G$ and
load impedance $Z_L$; \textup{(b)} circuit of Figure~\protect\ref{ports-fig}\textup{(a)}
with a two-port inserted between generator and load, depicting
input ($a_\cdot$) and output waves ($b_\cdot$) at the two-port
and at load.}\label{ports-fig}
\end{figure}

We begin our discussion with a \textit{one-port}, a single load
impedance serially inserted in a circuit of a voltage generator and its
impedance [see Figure \ref{ports-fig}(a)].
%
%
%
Recall that voltages, currents and impedances in an alternating current
(AC) circuit are modeled by complex numbers; otherwise, a loss of
alternating real voltage $u\dvtx t\mapsto u_0 \cos(\omega t)$ over a load
giving a phase shifted current $j\dvtx t \mapsto j_0 \cos(\omega t +\phi)$
would result into an awkward time dependent real resistance
$u(t)/j(t)$. In complex notation, the ratio of voltage $u(t)= u_0
e^{i\omega t}$ over current $j(t)= j_0 e^{i(\omega t+\phi)}$ is
constant and called \textit{impedance}
\[
Z:= \frac{u(t)}{j(t)} = \frac{u_0}{j_0} e^{i\phi} \in\mathbb C .
\]
Its real part is called \textit{resistance}, the imaginary part is the
\textit{reactance}. For example, under an AC-voltage $u_0\cos(\omega
)t$, a
serial circuit of a resistor with direct current DC-resistance $R$ and
an ideal capacitor with capacitance $C$ features an inverse impedance
(called \textit{admittance}) of $Z^{-1} = R^{-1}+i\omega C$. In fact, in
realistic scenarios, the resistance is positive, thus $Z \in-i\mathbb
H$, and the boundary (the imaginary axis) corresponds to ideal
(lossless) impedances.

\subsection{Reflections and characteristic impedance}\label{Reflection_scn}
We now assume that our circuit features a generator generating the
\textit{open circuit voltage} $u_G$ with internal impedance $Z_G$ and a
load with impedance $Z_L$ as depicted in
Figure \ref{ports-fig}(a)
with total impedance $Z= Z_G+Z_L$ according to Kirchhoff's circuit law.

Inspired by the wave model, voltage loss $u_L$ over and current flow
$j_L$ along the load is considered to be the superimposition of an
incoming (denoted by ``$+$'') and a reflected wave (denoted by ``$-$'')
in such a way that each single wave satisfies Ohm's law with a common
characteristic impedance $Z_c$. Since the reflected wave propagates
into a direction opposite to the incoming wave, we have the \textit{ansatz}
\[
u_L = u_L^+ + u_L^- \quad\mbox{and}\quad j_L = j^+_L - j^-_L
\]
with Ohm's law
\[
\frac{u^+_L}{j^+_L} = Z_c = \frac{u^-_L}{j^-_L} \quad\mbox{and}\quad \frac
{u^+_L + u^-_L}{j^+_L-j^-_L} = Z_L .
\]
The specific decomposition or equivalently the choice of $Z_c$ is
arbitrary in many applications, and it will be guided by imposing
additional conditions. Usually $Z_c$ is taken positive, or at least
chosen such that the normalized impedances $\widetilde{Z} := Z/Z_c$
will be again of positive real part, that is, $\widetilde{Z} \in
-i\mathbb H$. The analog for the right half plane of the Cayley
transform, (\ref{Cayley-transform}), yields then \textit{reflection
coefficient} of the load
%
%
\begin{equation}\label{load-refl}
\Gamma_L := \frac{u^-_L}{u^+_L} = \frac{j^-_L}{j^+_L} = \frac
{Z_L-Z_c}{Z_L+Z_c} = \frac{\widetilde{Z}_L-1}{\widetilde
{Z}_L+1}=\mathcal{C}(-i\widetilde{Z}_L)
\end{equation}
as an element of the unit-disk $\mathbb D$. Of course, there is no
reflection if $Z_L=Z_c$.

\subsection{The chain matrix}
We are now in a position to investigate the generic scenario of signal
transmission through a \textit{two-port} [see Figure \ref{ports-fig}(b)].
Among others due to linearity of the Maxwell equations, voltages and
currents $(j_1=a_1,j_2=a_2)$ on either side of the two-port have a
linear relationship governed by a so-called impedance matrix $Z$
\[
Z\pmatrix{j_1\cr j_2}
=\pmatrix{Z_{11}&Z_{12}\cr Z_{21}&Z_{22}}
\pmatrix{j_1\cr j_2}
= \pmatrix{u_1\cr u_2}.
\]
For given circuit parameters, the coefficients of the impedance matrix
can be easily computed. For example, $Z_{11} = \frac{u_1}{j_1}
|_{j_2=0}$ is the well-known input-impedance (by inserting a load of
infinite impedance the right-hand side becomes an open circuit with
$j_2=0$). In most applications it turns out that $Z$ is symmetric; the
corresponding two-port is then called \textit{reciprocal}.

One easily verifies that the \textit{chain matrix} which is usually
denoted by \cha\ (the Russian letter ``cha'') relating $(u_2,-j_2)$
with $(u_1,j_1)$ is given by
\[
\mbox{\cha}\pmatrix{u_2\cr-j_2}
=\frac{1}{Z_{21}}\pmatrix{Z_{11}&\det(Z)\cr1& Z_{22}}
\pmatrix{u_2\cr-j_2}
=\pmatrix{u_1\cr j_1}
\]
[in contrast to the mathematical literature, the engineering literature
tends to use a \textit{transmission matrix} relating $(u_2,j_2)$ with
$(u_1,j_1)$ instead, with reversed $j_2=-a_2$ in Figure \ref{ports-fig}(b)].
An advantage of the chain matrix over the impedance matrix is that the
former is well defined for the limit $Z_{21}\to\infty$, for example,
for an ideal coil in series with the load.

Again, only using lossless (i.e., purely imaginary) impedances (such
as ideal inductances, transformers and capacitors) guarantees that the
corresponding chain matrix has real diagonal coefficients and imaginary
coefficients elsewhere. Moreover, for cascaded two-ports (i.e.,
several two-ports in serial connection), the resulting chain matrix is
just the product of the individual chain matrices. By linear algebraic
decomposition of $\mathbb{S}\mathbb{L}(2,\mathbb R)$ it can be shown
that every
lossless two-port of lumped elements can be modeled by cascading
combinations of two-ports involving only inductances, transformers and
capacitors (see \cite{H82}, page 18).
Note that $j_2 = a_2=-a_L=-j_L$ in Figure \ref{ports-fig}(b).
Hence \cha\ defines a M\"obius transformation relating load impedance
with the input impedance of the two-port
\[
M_{\mbox{\chaa}}(Z_L) = \frac{\mbox{\cha}_{11}Z_L + \mbox{\cha
}_{12}}{\mbox{\cha}_{21}Z_L +
\mbox{\cha}
_{22}} = \frac{\mbox{\cha}_{11}u_2 - \mbox{\cha}_{12}j_2}{\mbox
{\cha}_{21}u_2 -
\mbox{\cha}
_{22}j_2} = \frac{u_1}{j_1} = Z_1 .
\]
Here, $Z_1 = u_1/j_1$ is the impedance of the load $Z_L=-u_2/j_2$ as
viewed through the two-port.
As a consequence we make the following remark.
\begin{Rm} Serial cascading of lossless two-ports is equivalent to the
action of the M\"obius group $\mathbb{S}\mathbb{L}(2,\mathbb R)$ on
the $i$-fold
$iZ_L\in\mathbb H$ of load impedances $Z_L$.
\end{Rm}

We are thus led to the statistical inverse problem (cf. Problem \ref
{prob_statement}).
\begin{Prob}\label{est-imp-prob}
Estimate the load impedance $Z_L$ when only the impedance $Z_{(1)}$
viewed through the two-port
can be observed
where $
\varepsilon= \mathcal{I}\circ M_{\mbox{\chaa}}\circ\mathcal
{I}^{-1} \in\mathbb{S}\mathbb{L}
(2)$ is assumed to be of known distribution. Here, $\mathcal{I}
\dvtx\mathbb C\to\mathbb C\dvtx z \to iz$ denotes the multiplication with $i$.
\end{Prob}

In conclusion we note that one may as well consider normalized
impedances $\widetilde{Z}=Z/Z_c$ or equivalently reflection
coefficients, (\ref{load-refl}). Then the mapping for the normalized
impedances goes as follows:
\[
\frac{Z_1}{Z_c} = \frac{\mbox{\cha}_{11} {Z_L}/{Z_c} +
{\mbox{\cha}
_{12}}/{Z_c} }{Z_c\mbox{\cha}_{12} {Z_L}/{Z_c} + \mbox{\cha
}_{22}} .
\]

\subsection{Estimating resistances seen through electrolyte capacitors}
\label{example-scn}

It is well known that over the duration of years properties of
electronic equipment change due to wear-out effects of various
elements. In particular, electrolyte capacitors have a tendency to
loose capacitance. In effect, older electronic devices deviate from
original calibration and may feature nondesired side-effects; for
example, field strengths of transmitters may grow stronger than tolerated.

%
%
\begin{figure}

\includegraphics{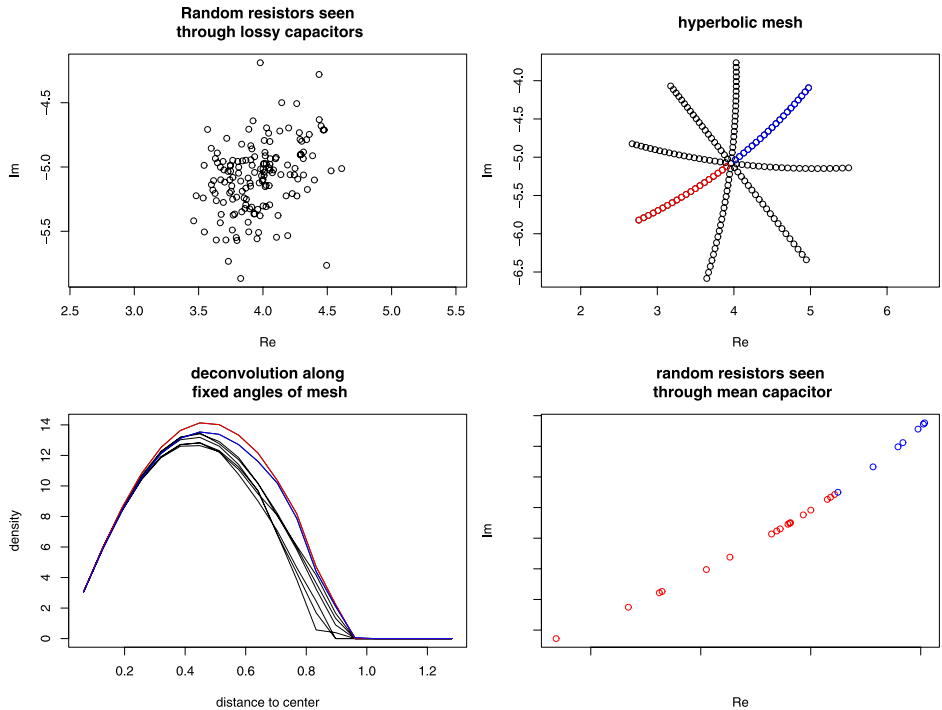}

\caption{Random resistors seen through random
capacitors. Left top: original measurements. Right top: hyperbolic
polar mesh points at which the M\"{o}bius deconvolution was computed.
The bottom left image depicts the deconvolved densities times
hyperbolic measure along fixed angles. The red and blue line goes along
the red and blue mesh points. For verification, the one-dimensional
distribution of the resistors seen through the single mean capacitor is
depicted in the bottom right image.}
\label{resistor_capacitor_20rads_fig}
\end{figure}

In an application of our method we consider a series of $n=150$
measurements of random resistors of $15$ $\Omega$ provided with an
accuracy of 10 percent by the manufacturer (they range from $13.5$
$\Omega$ to $17.7$ $\Omega$) viewed through 30 random capacitors at 1
kHz taken at the Department of Electrical Engineering, University of
Applied Sciences, Fulda, Germany. These originally identical lossy
$22$ $\mu $F capacitors have been collected from over ten year old
electronic gear. For the impedance measurements the LCR-Bridge
``HM8118'' has been used that comes with an accuracy of $0.3 \%$
guaranteed by its producer HAMEG. We model the $i$-fold of the
impedance $Z$ of these capacitors with a hyperbolic
Gaussian-distribution at unit impedance, that is, $Z i=Z_cM_{\rho
_{\varepsilon}}(i)$ with a suitable characteristic impedance $Z_c$ and
a random hyperbolic Gaussian $\rho$-distributed M\"obius transformation
$M_{\rho_\varepsilon}$. Measurement of the capacitors gives $\rho
_\varepsilon\approx0.0004$ corresponding to a spread of roughly $4.8
\%$. Our goal lies in the
reconstruction of the one-dimensional resistances $R$ solely from the
observations
\[
W= \frac{1}{{1}/{R} + {1}/{Z}} = \frac{ZR}{Z+R}
\]
and the known dispersion $\rho_\varepsilon$ of the corruption as posed
in Problem \ref{est-imp-prob}. To this end we apply M\"{o}bius
deconvolution to the model
\[
\frac{W}{W_c} i = M_{\varepsilon} \biggl(\frac{R}{R_c} i \biggr)
\]
with suitable characteristic impedances $W_c$ and $R_c$. For this
application the Euclidean means have been chosen as characteristic
impedances. Alternatively, a better approach may be to
use hyperbolic intrinsic means (see \cite{BP03} and \cite{BP05}).
Figure \ref{resistor_capacitor_20rads_fig}
shows the observations $W$ in the left top corner. M\"{o}bius
deconvolution is computed along the hyperbolic polar mesh-points
depicted in the right top corner
(Figure~\ref{resistor_capacitor_20rads_fig}).
Below in the bottom left corner the deconvolved densities along
fixed angles of the mesh are depicted. The angle depicted in red shows
highest density followed by the angle depicted in blue. The location of
the two dominating directions depicted with the same colors in the
right top corner
(Figure \ref{resistor_capacitor_20rads_fig})
is in high agreement with the location of the impedances of the
resistors seen through the mean capacitor depicted in the bottom right
image of Figure \ref{resistor_capacitor_20rads_fig}.
Indeed, one can say that with the few measurements available, we were
able to reconstruct the nature of the unobserved elements, namely
resistors with impedances distributed along a one-dimensional subset in
the complex half plane.

%
%
\begin{figure}

\includegraphics{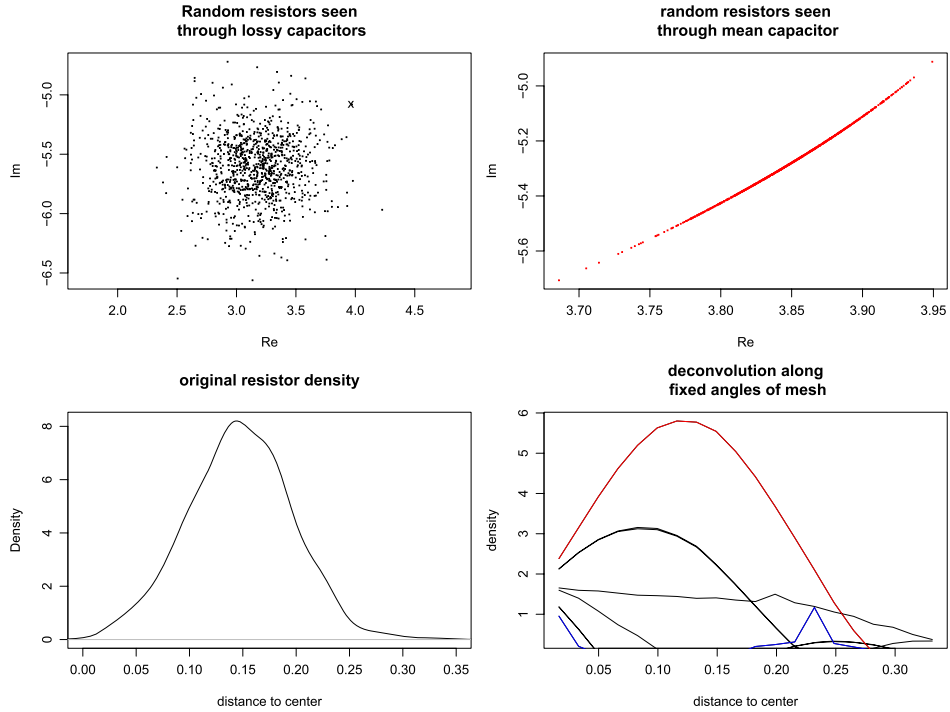}

\caption{Simulation of $n=1000$ random resistors seen through $n$
random capacitors. Top left: observed impedances, characteristic
impedance denoted by ``$x$.'' Top right: unobserved underlying normalized
resistors seen through unobserved mean capacitor. Bottom left: radial
density times hyperbolic measure of nonobserved resistors. Bottom
right: deconvolved densities times hyperbolic measure along fixed
angles of mesh in Figure \protect\ref{resistor_capacitor_20rads_fig}.}
\label{simu_resistor_capacitor_20rads_fig}
\end{figure}

In Figure \ref{simu_resistor_capacitor_20rads_fig}
the above scenario is more prominently reenacted in a simulation of
$n=1000$ measurements using $R_c=\min(R)$ and $W_c$ from the preceding
example (depicted by ``$x$''). We show observed measurements and
unobserved capacitors in the top row as well as original density and
the deconvolved densities along respective angles as in Figure \ref
{resistor_capacitor_20rads_fig}
in the bottom row. Obviously the distribution of the unobserved
resistors is quite reasonably recovered along the grid.

\begin{appendix}

\section{Polar coordinates and convolution}
\label{upper-half-plane-scn}

In this appendix we
focus on the right and left action (they are different due to
noncommutativity) of the special orthogonal group $\mathbb{S}\mathbb
{O}(2)$ on $\mathbb{S}\mathbb{L}(2)$
giving rise to polar coordinates and to
$\mathbb H$ viewed as the quotient with respect to one of the actions.
In fact, the action of $\mathbb{S}\mathbb{L}(2)$ on $\mathbb H$ can
be naturally viewed
in polar coordinates (cf. Figure \ref{hyperbolic:fig}) of which we will
make extensive use in the proof of Theorem \ref{thm:LowerBound} in
Appendix \ref{ssec:LowerBound}.

%
%
\begin{figure}

\includegraphics{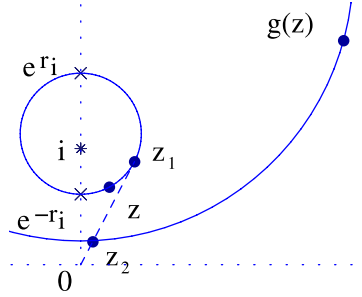}

\caption{The polar coordinates $(r,u)$ of $z\in\mathbb H$
are obtained from the unique circle hyperbolically centered at $i$
(star) containing $z$; that is, this circle is symmetric to the
imaginary axis and intersects it at points of form $e^{-r}i$ and
$e^{r}i$ ($x$-crosses). Rotating $z$ by the hyperbolic angle $-u$ along
this circle, the lower point $e^{-r}i$ is obtained. The polar\vspace*{1pt}
coordinates $(u_2,R,u_1)$ of $g\in\mathbb{S}\mathbb{L}(2)$ rotate
$z$ by the hyperbolic
angle $u_1$ along the above circle to obtain $z_1$, rescale $z_1$ by
$e^{-R}$ to obtain $z_2=e^{-R}z_1$ and subsequently rotate $z_2$ by the
hyperbolic angle $u_2$ along the unique hyperbolic circle through $z_2$
hyperbolically centered at $i$ yielding $g(z)$.}\label{hyperbolic:fig}
\end{figure}

In contrast to a metric on a manifold in the usual topological
sense, a \textit{Riemannian metric} is a metric in every tangent
space, that varies smoothly with the offset of the tangent
space. Thereby, every Riemannian metric defines a metric on the
manifold in the usual sense and a unique volume element, called
the \textit{Riemannian volume}, giving rise to a unique measure
on the manifold. For more details we refer to Lee
\cite{Lee97}, Chapter 3.

If the Riemannian metric on $\mathbb{S}\mathbb{L}(2)$ underlying the
Haar measure $dg$
is chosen such that the natural Riemannian quotient metric on $\mathbb
H$ yields the hyperbolic measure $dz$ on $\mathbb H$, then we are able
to lift densities on the bottom space $\mathbb H$ to the top space
$\mathbb{S}\mathbb{L}
(2)$ to obtain (\ref{conv_lift:def}) yielding (\ref{conv_model}).

Let us begin with the observation that the hyperbolic measure (\ref
{hyperbolic-measure}) is the Riemannian volume element
of the extension of the standard Euclidean metric in the tangent space
of $\mathbb H$ at $z=i$ in a left-invariant way under the action of
$\mathbb{S}\mathbb{L}
(2)$ on~$\mathbb H$.
Similarly, we equip $\mathbb{S}\mathbb{L}(2)$ with the Riemannian
metric obtained from
the left $\mathbb{S}\mathbb{L}(2)$-invariant extension of the
standard Euclidean metric
in the tangent space of the unit matrix $\mathbf{I}\in\mathbb{S}\mathbb{L}(2)$.
We denote the corresponding Riemannian volume element which defines a
left-invariant Haar measure by $dg$.
According to \cite{T85}, Exercise 19, page 149, it is also right
invariant, that is,
\[
\int_{\mathbb{S}\mathbb{L}(2)} f(agb) \,dg =\int
_{\mathbb{S}\mathbb{L}(2)} f(g) \,dg\qquad
\forall a,b\in\mathbb{S}\mathbb{L}(2),f \in L^1 (\mathbb
{S}\mathbb{L}(2),dg ) .
\]
As mentioned before, for arbitrary $z,z'\in\mathbb H$ there exists a $g
\in\mathbb{S}\mathbb{L}(2)$ such that $z' = M_g (z)$. Given one such
$g$, any other $g'
\in\mathbb{S}\mathbb{L}(2)$ satisfies $M_{g'} (z) = z'$ if, and only
if, $g^{-1} g'
\in
\mathbb{S}\mathbb{O}(2)$. In particular, $M_g (i) = i$ if, and only
if, $g \in\mathbb{S}\mathbb{O}(2)$.
This entails that the following mapping of the quotient space $\mathbb
{S}\mathbb{L}
(2)/\mathbb{S}\mathbb{O}(2)$ due to the right action of $\mathbb
{S}\mathbb{O}(2)$ is well defined and bijective
%
%
\begin{equation}\label{SL2-H-isometry}
\left.
\begin{array}{rcl} \mathbb{S}\mathbb{L}(2)/\mathbb{S}\mathbb{O}(2)
&\to& \mathbb H\\
g\mathbb{S}\mathbb{O}(2) &\mapsto& M_g(i)
\end{array}
\right\} .
\end{equation}
Since the mapping preserves the action of $\mathbb{S}\mathbb{L}(2)$,
the natural
Riemannian quotient metric of $ \mathbb{S}\mathbb{L}(2)/\mathbb
{S}\mathbb{O}(2)$ is isometric with the
the hyperbolic metric of $\mathbb H$.

Next, consider the left action of $\mathbb{S}\mathbb{O}(2)$ on
$\mathbb{S}\mathbb{L}(2)$. This projects
to a left-action on $\mathbb H$ giving rise to {polar coordinates} $u
\in[0,2\pi)$, called the \textit{hyperbolic angle} and $r>0$ (cf. Figure
\ref{hyperbolic:fig}) of
\[
z = M_{k_{u}}\circ M_{R_{r}} (i)\qquad\mbox{with } k_{u} = \pmatrix{
\cos u& \sin u\cr-\sin u &\cos u},
R_{r}=\pmatrix{e^{-r/2}&0\cr0&e^{r/2}}.
\]
For ease of notation, for the entire paper we identify $M_{k_{u}}$ with
$k_u$ and $M_{R_{r}}$ with $R_r$, respectively, such that
$k_u \circ R_r (i) = k_u (e^{-r}i)$.
Since for every $z\in\mathbb H\setminus\{i\}$ we have that $k_{\phi
}(z) = z$ if, and only if $\phi\equiv0\operatorname{mod}\pi$,
polar coordinates
cover the hyperbolic plane twice. As a consequence, $z\in\mathbb H$ %
has polar coordinates with $r\geq0$ uniquely determined and $u$
unique modulo $\pi$ if $r>0$. Thus, the hyperbolic area element (\ref
{hyperbolic-measure}) transforms to
%
%
\begin{equation}\label{hyperbolic-measure-polar-coordinates}
dz = \sinh(r) \,dr \,du .
\end{equation}
Polar coordinates can also be defined on $\mathbb{S}\mathbb{L}(2)$:
every element $g
\in
\mathbb{S}\mathbb{L}(2)$ has a decomposition
%
%
\begin{equation}\label{SL(2)-polar-decomposition}
g = k_u R_r k_{u'},\qquad r \in[0,\infty), u,u' \in[0,2\pi),
\end{equation}
with uniquely determined $r\geq0$; if $r>0$ then $u$ and $u'$ are also
uniquely determined modulo $\pi$ (in fact, one of the two is unique
modulo $2\pi$). In view of the isometry (\ref{SL2-H-isometry}), this
gives our choice of Haar measure in polar coordinates
\[
dg = \sinh(r)\,
dr \,du \,du' .
\]

If $g, g'$ are independent random elements in $\mathbb{S}\mathbb
{L}(2)$ with densities
$f_1, f_2$ continuous with
respect to Haar measure, we have for the probability that the product
is contained in a measurable subset $A\subset\mathbb{S}\mathbb
{L}(2)$ by
left-invariance of the measure that $\mathbb{P}(gg'\in A) = \int_{A}
(f_1*f_2)(a) \,da$ with
the convolution of $f_1$ and $f_2$ given by
\[
(f_1*f_2)(a) := \int_{\mathbb{S}\mathbb{L}(2)} f_1(g)f_2(g^{-1}a)
\,dg .
\]
In general, convolutions over noncommutative groups are noncommutative.

Suppose now that $Z$ is a random quantity on $\mathbb H$ with density
$f_2$ continuous with respect to the hyperbolic measure. Using polar
coordinates, this density lifts to a right $\mathbb{S}\mathbb
{O}(2)$-invariant density
$\widetilde{f}_2$ on $\mathbb{S}\mathbb{L}(2)$
$\widetilde{f}_2(k_u R_r k_u') := f_2 (k_uR_r(i) )$.
Hence, convolutions of a density $f_1$ on $\mathbb\mathbb{S}\mathbb
{L}(2)$ with a
density $f_2$ on $\mathbb H$ can be well defined by lifting to a right
$\mathbb{S}\mathbb{O}(2)$-invariant density on $\mathbb{S}\mathbb{L}(2)$
%
%
\begin{equation}\label{conv_lift:def}
(f_1*f_2)(z) := \int_{\mathbb{S}\mathbb{L}(2)} f_1(g) \widetilde
{f}_2(g^{-1} a) \,dg
\end{equation}
with any $a\in\mathbb{S}\mathbb{L}(2)$ giving $M_a(i) = z$.
This convolution is commutative if either $f_{1}$ is bi-invariant or
if $f_2$ is $\mathbb{S}\mathbb{O}(2)$-invariant.

\section{Proofs}\label{sec:proofs}

\subsection{\texorpdfstring{Upper bound: Proof of Theorem
\protect\ref{main-thm}}{Upper bound: Proof of Theorem 3.1}}
\label{ssec:upper}
In order to measure the performance of $f^{(n,T)}_{X}$ we consider the
mean integrated squared error
\[
\mathbb{E}\bigl\|f_X-f_{X}^{(n,T)} \bigr\|^2 = \mathbb{E}\bigl\| f_X - \mathbb{E}
f_{X}^{(n,T)} \bigr\|^2+ \bigl\|\mathbb{E}f_{X}^{(n,T)} - f_X \bigr\|^2
\]
with the usual variance-bias decomposition. The assertion of Theorem
\ref{main-thm} then follows from the following more detailed lemma.
%
%
\begin{Lem}\label{var_bias:lem} For $f_X, f_Y$ and $f_{\varepsilon}$
satisfying \textup{(D.1)} and \textup{(D.2)}, and $\mathcal{
H}f_{\varepsilon
}$ bounded from below on compact sets,
there is a constant $C>0$ independent of $T$ and $n$ such that
\[
\mathbb{E}\bigl\| f_X - \mathbb{E}f_{X}^{(n,T)} \bigr\|^2 \leq\frac{C}{\inf_{|t|<T}
{|\mathcal{H}f_{\varepsilon}({1}/{2} +it)|^2}}\frac{T^2}{n} .
\]
If, additionally, $f_X$ satisfies \textup{(D.3)}, then
\[
\bigl\|\mathbb{E}f_{X}^{(n,T)} - f_X \bigr\|^2 \leq QT^{-2\alpha} .
\]
\end{Lem}
\begin{pf} We first note that by (\ref{inv_transfo}), definition
(\ref{inv_emp_transfo}), since the right-hand side is in $L^2(\mathbb H,dz)$,
\[
\mathcal{H}f_X^{(n,T)} \biggl(\frac{1}{2}+it,k \biggr) = \frac{\mathcal
{H}f_Y^n({1}/{2}+it,k)}{\mathcal{ H}f_{\varepsilon}({1}/{2}+it,k)} {\mathbb I}_{(-T,T)}(t) .
\]
Here $\mathbb I$ denotes the \textit{indicator function}.
Hence by the Fubini--Tonelli theorem, (\ref{inv_emp_transfo}),
(\ref{emp_transfo_exp_value}), (\ref{inv_transfo}) and (\ref
{helg_trans_conv})
%
%
\begin{eqnarray}\label{HEf_X_nT}\quad
&&
\mathcal{H} \bigl(\mathbb{E}f_X^{(n,T)}\bigr) \biggl(\frac{1}{2}+it,k \biggr)\nonumber\\
&&\qquad= \mathcal{H} \biggl(z\to\int_{|t'|< T}\int_{u'=0}^{u'=2\pi
} \frac{\mathbb{E}\mathcal{ H}f^{n}_Y ({1}/{2} +it',k_{u'}
)}{\mathcal{ H}f_{\varepsilon} ({1}/{2} +it' )} \nonumber\\
&&\qquad\quad\hspace*{100.8pt}{}\times(\operatorname{Im}
(k_{u'}(z) ) )^{{1}/{2}+it'}
\,d\tau' \biggr) \biggl(\frac{1}{2}+it,k \biggr) \nonumber\\
&&\qquad=\mathcal{H}f_X \biggl(\frac{1}{2}+it,k \biggr) {\mathbb
I}_{(-T,T)}(t) .
\end{eqnarray}
Deduce from (\ref{emp_transfo}),
%
%
\begin{eqnarray}\label{EHf_Yn}
&&\mathbb{E}\biggl|\mathcal{H}
f^{(n)}_Y \biggl(\frac{1}{2}+it,k \biggr)
\biggr|^2\nonumber\\
&&\qquad
= \biggl|\mathcal{H}f_Y \biggl({\frac{1}{2}+it},k \biggr) \biggr|^2 +
\frac{\operatorname{Im}(\mathbb{E}k(Y) )- |\mathcal{H}f_Y ({
{1}/{2}+it},k ) |^2}{n} .
\end{eqnarray}
In addition, $\operatorname{Im}(\mathbb{E}k(Y) ) = \mathcal
{H}f_Y(1,k)$ implies that
%
%
\begin{equation} \label{est_int_Ek}
{\int_0^{2\pi}} |{\operatorname{Im}}(\mathbb{E}k_u(Y) ) | \,du =
{\int_0^{2\pi}} |\mathcal{H}f_{Y}(1,k_u)| \,du\leq C
\end{equation}
with a suitable constant $C>0$, since $\mathcal{H}f_Y$ is analytic.
Thus, using the Plancherel identity (\ref{planch_eq}),
(\ref{helg_trans_conv}), the Fubini--Tonelli
theorem, (\ref{emp_transfo_exp_value}),
(\ref{HEf_X_nT}), (\ref{EHf_Yn}) and (\ref{est_int_Ek}) [by hypothesis
(D.1) and
(\ref{planch_eq}), $\|f_Y\|^2 = \|\mathcal{H}f_Y\|^2 < \infty$], we
have indeed for the variance
\begin{eqnarray*}
&& \mathbb{E}\int_{\mathbb H} \bigl|f_X^{(n,T)} - \mathbb
{E}f^{(n,T)}_X\bigr|^2 \,dz\\
&&\qquad=\int_{t<|T|}\int_{u\in[0,2\pi)} \frac{\mathbb{E}|\mathcal
{H}f^{(n)}_Y -
\mathcal{H}f_Y|^2}{|\mathcal{H}f_{\varepsilon}|^2} \,d\tau\\
&&\qquad=\frac{1}{n} \frac{1}{8\pi^2} \int_{t<|T|}\int_{u\in[0,2\pi)}
\frac{\operatorname{Im}(\mathbb{E}k_u(Y)) - |\mathcal{H}f_Y({1}/{2}
+it)|^2}{|\mathcal
{H}f_{\varepsilon}({1}/{2} +it)|^2}
t \tanh(\pi
t) \,dt \,du\\
&&\qquad\leq \frac{C}{\inf_{|t|<T}
{|\mathcal{H}f_{\varepsilon}({1}/{2} +it)|^2}} \frac{T^2}{n}
\end{eqnarray*}
with a constant $C>0$ involving
neither $n$ nor $T$.

In the next step we similarly estimate the squared bias under the
additional assumption
(D.3)
using also (\ref{planch_eq}) and (\ref{HEf_X_nT})
\begin{eqnarray*}
&& \int_{\mathbb H} \bigl|\mathbb{E}f_X^{(n,T)} - f_X\bigr|^2 \,dz
\\
&&\qquad=\int_{|t|\geq T}\int_{u\in[0,2\pi)} |\mathcal{ H}f_X|^2 \,d\tau\\
&&\qquad= \frac{1}{8\pi^2}\int_{|t|\geq T}\int_{u=0}^{u=2\pi}
\biggl| \biggl(\frac{1}{2}+it \biggr) \biggl(-\frac{1}{2}+it \biggr) \biggr|^{-\alpha}
\biggl| \biggl(\frac{1}{2}+it \biggr) \biggl(-\frac{1}{2}+it \biggr)\biggr|^{\alpha}\\
&&
\hspace*{112.3pt}{}\times\biggl|\mathcal{ H}f_X \biggl(\frac{1}{2} +it,k_u \biggr) \biggr|^2 t
\tanh(\pi t) \,dt \,du
\\
&&\qquad\leq QT^{-2\alpha} .
\end{eqnarray*}
\upqed\end{pf}

\subsection{\texorpdfstring{Optimal rate: Proof of Theorem
\protect\ref{opt-T-cor}}{Optimal rate: Proof of Theorem 3.2}}

If $f_{\varepsilon}$ satisfies (D.4) we have the upper bound
%
%
\begin{equation}\label{var:bound:log}
\log C + \log C_1+ 2\log T -\log n + \frac{2}{\gamma} T^\beta
\biggl(1+\frac{1}{4T^2} \biggr)^{{\beta}/{2}}
\end{equation}
for the logarithm of the variance term (cf. Lemma \ref{var_bias:lem}).
A sufficient condition for convergence of the variance term while
$T=T(n) \to\infty$, is that (\ref{var:bound:log}) tends to $-\infty$.
Hence, $T$ is of form
\[
T(n) = \biggl(\frac{\gamma}{2} \log n - \frac{\gamma}{2} A(n)
\biggr)^{1/\beta}
\]
with $A(n) \to+\infty$ at a rate lower than that of $\log n$. A short
computation gives the rate
\[
\frac{2}{\beta} \log(\log n)-A(n)
\]
for (\ref{var:bound:log}). In case of optimality this rate must be
larger or equal to the logarithmic rate of the upper bound of the bias
term in Lemma \ref{var_bias:lem} which is then
\[
-\frac{2\alpha}{\beta} \log(\log n) .
\]
In consequence the rate of $A(n)$ is $\eta\log(\log n)$ with $\eta
\geq
2(1+\alpha)/\beta$ as asserted by Theorem \ref{opt-T-cor}.

\subsection{\texorpdfstring{Lower bound properties: Proof of Theorem
\protect\ref{thm:LowerBound}}{Lower bound properties: Proof of Theorem 3.3}}
\label{ssec:LowerBound}

Recall the decomposition in polar coordinates from Appendix \ref
{upper-half-plane-scn} (cf. Figure \ref{hyperbolic:fig}).
The general idea of proof goes as follows.
Define the \textit{dilation} of $H$ by
%
%
\begin{equation}\label{eq:Hdelta}
H^{\delta} (k (e^{-r} i) ) :=
H(k( e^{-\delta r}i)) \frac{\sinh(\delta r)}{\sinh r}
P_{-{1}/{2}} (\cosh(\delta r) ).
\end{equation}
For arbitrary $\mathbb{S}\mathbb{O}(2)$-invariant $H$ and $\delta
>0$, $H^\delta$ is
obviously also $\mathbb{S}\mathbb{O}(2)$-invariant

To derive a lower bound for estimating $f_X$ in $L^2$ norm, we follow a
classical scheme which has been condensed in \cite{Fan91}, pages
1261 and 1262 (cf. also \cite{HM07}, pages 1555 and 1556). The
adaption of this scheme to the Poincar\'e plane, however, is not at all
obvious and will be the subject of the following sections. After some
elaborate preparation in the following two sections we take a pair $f_0
\in{ {\mathcal{F}_{\alpha}(Q)}}$, $f_n \in{ {\mathcal{F}_{\alpha
}(Q)}}$, for which
\[
f_n = f_0 + C_H \delta^{-\alpha} H^\delta,
\]
where $\delta= \delta_n$ (cf. Sections \ref{f0:scn} and \ref
{feps:scn} below).
Then, in Section \ref{chi2:dist:scn} we show that $\delta$ can be
chosen such that
%
%
\begin{equation}\label{eq:chisq-1}\quad
\chi^2(f_\varepsilon* f_0 , f_\varepsilon* f_n) := \int_0^\infty
(f_\varepsilon* f_0 - f_\varepsilon* f_n)^2 (f_\varepsilon*
f_0)^{-1} \,dz \le\frac C n .
\end{equation}
In consequence by (3.3) of \cite{Fan91}, there is $d_1>0$ such that for
any estimator $f^n$ of $f_X$,
%
%
\begin{equation}\label{eq:LowerP}
\sup_{f_X \in\{f_0, f_n\}} {\mathbb P}_f \{ \| f^n - f_X \|_2 > \|
f_0 - f_n \|_2/2 \} > d_1,
\end{equation}
which gives with (3.4) of \cite{Fan91} a lower bound
\[
\sup_{f_X \in{\mathcal{F}_{\alpha}(Q)}}E\|f^n-f_X\|^2 \geq\frac
{d_1}{4} \| f_0 - f_n \|_2
\asymp\delta^{-\alpha} \| H^\delta\|_2.
\]
The choice of $\delta$ in (\ref{eq:delta-rate}) at the end of Section
\ref{chi2:dist:scn} in conjunction with $\|H^\delta\|^2\asymp\|H\|^2$
from Lemma \ref{lem:Dilate} then yields the rate
\[
\delta^{-\alpha} \| H^\delta\|^2\asymp(\log n )^{-2
{\alpha}/{\beta}}
\]
asserted by Theorem \ref{thm:LowerBound}.

\subsubsection{Convolution equation in polar coordinates}

\begin{Lem}\label{lem:conv.rep}
Suppose $f$ is bi-invariant and $h$ is $\mathbb{S}\mathbb{O}(2)$-invariant.
Write $z = k_u R_r(i)$, $g = k_\alpha R_s k_\beta\in\mathbb
{S}\mathbb{L}(2)$
and $dg =
\sinh s\, d\alpha \,ds \,d\beta$.
Then,
\[
(f * h)(z) = 2\pi\int_{\phi=0}^{2\pi} \int_{s=0}^\infty
f(e^{-s} i) h (e^s k_{-\phi} (e^{-r} i) ) \sinh s \,d\phi \,ds.
\]
\end{Lem}
\begin{pf}
Since $k_{-\alpha}k_u = k_{u-\alpha}$ observe that
\begin{eqnarray*}
(f * h)(z) &=& \int_{\mathbb{S}\mathbb{L}(2)} f(g) h(g^{-1} z) \,dg \\
&=& \dint_{\alpha,\beta\in[0,2\pi)}\int_{s=0}^{s=\infty}
f(k_\alpha R_s
k_\beta)
h (k_{-\beta} R_{-s} k_{-\alpha} k_u R_r(i) )\\
&&\hspace*{84.6pt}{}\times
\sinh s \,d\alpha \,ds \,d\beta\\
&=& 2\pi\int_{\phi=0}^{2\pi} \int_{s=0}^\infty f(e^{-s} i)
h (R_{-s} k_{-\phi} R_r(i) ) \sinh s \,d\phi\, ds .
\end{eqnarray*}
\upqed\end{pf}
\begin{Lem}\label{lem:R}
Define $\eta(r,s,\phi)$ and $R(r,s,\phi)$ by $k_{\eta(r,s,\phi)}
e^{-R(r,s,\phi)}i = e^s k_{-\phi} \times\break e^{-r} i$,
where $0 \le\eta(r,s,\phi) < 2\pi$ and $R(r,s,\phi) \ge0$.
Suppose $\phi\in[0,2\pi)$ and $r, s \ge0$.
Then,
\[
|r-s| \le R(r,s,\phi) \le r+s.
\]
\end{Lem}
\begin{pf}
Let $\psi= 2\phi$.
From \cite{T85}, page 125, we take
\[
k_{-\phi} e^{-r} i =
\frac{-\sin\psi\sinh r + i}{\cos\psi\sinh r + \cosh r} .
\]
Let $k_{\eta(r,s,\phi)} e^{-R(r,s,\phi)}i = x +i y$.
Then using \cite{T85}, page 150,
\begin{eqnarray*}
\cosh( R(r,s,\phi) )
&=& \frac{1}{2} \frac{1+x^2+y^2}{y} \\
&=& \frac{1}{2} \biggl(e^{-s}(\cos\psi\sinh r + \cosh r) +
e^s \frac{\sin^2\psi\sinh^2 r +1}{\cos\psi\sinh r + \cosh
r} \biggr).
\end{eqnarray*}
Set $t = \cos\psi$. Since
\begin{eqnarray*}
\frac{\sin^2 \psi\sinh^2 r + 1}{\cos\psi\sinh r + \cosh r}
&=& \frac{(1 - t^2) \sinh^2 r + 1} {t\sinh r + \cosh r } \\
&=& \cosh r - t \sinh r,
\end{eqnarray*}
we have
\begin{eqnarray*}
\cosh( R(r,s,\phi) )
&=& \tfrac1 2 \{ e^{-s} (\cosh r + t\sinh r) + e^s (\cosh r - t
\sinh r) \} \\
&=& -t \sinh r \sinh s + \cosh r \cosh s.
\end{eqnarray*}
Since $-1 \le t \le1$,
$\cosh( R(r,s,\phi) )$ has the maximum $\cosh(r+s) $ at $t =
-1$ and
the minimum $\cosh(r - s)$ at $t = 1$.
Suppose that $\cosh a \le\cosh b$ with $b \ge0$.
This implies that $b \ge|a|$ since
\[
\cases{
b \ge a, &\quad if $a \ge0$, \cr
b \ge-a, &\quad if $a < 0$.}
\]
The desired result follows from this and the assumption $R(r,s,\phi)
\ge0$.
\end{pf}

The following lemma is an immediate consequence of Lemma \ref{lem:R}.
\begin{Lem}\label{lem:b.m.H}
Suppose $H$ is $\mathbb{S}\mathbb{O}(2)$-invariant and $H(e^{-r}i)$
is monotonically
decreasing in $r$.
Then,
\[
H (e^s k_{-\phi} e^{-r} i ) \le H \bigl(e^{-|r-s|} i \bigr)
\qquad\mbox{for }s \ge0, \phi\in[0, 2\pi), r \ge0.
\]
\end{Lem}

\subsubsection{Dilation}

\begin{Lem} \label{lem:Dilate}
Suppose
%
%
\begin{equation}\label{eq:BoundedSupportofH}
\mathcal{H} H \bigl({\tfrac1 2}+ it \bigr) = 0 \qquad\mbox{for }t \notin\bigl[{\tfrac
1 2},
1 \bigr].
\end{equation}
Then for $\delta\to\infty$,
%
%
\begin{eqnarray}\label{eq:DilateH(t)}
\mathcal{H} {H^{\delta}} \biggl(\frac{1}{2}+it \biggr)
&\asymp&
\frac{1} {\delta} \mathcal{H} H \biggl(\frac{1}{2}+i \frac{t}{\delta
} \biggr),
\\
\|H^{\delta}\|^2 &\asymp& \|H\|^2,
\nonumber\\
\| \Delta^{\alpha/2} H^\delta\|_2
&\asymp& \delta^{\alpha/2} \| \Delta^{\alpha} H \|_2\nonumber
\end{eqnarray}
and for $f_{\varepsilon}$ satisfying \textup{(D.4)} with $\gamma=1$, with a
constant $C>0$,
\[
\int| f_\varepsilon* H^\delta|^2 \le C e^{-2 ({\delta
}/{2} )^{2\beta}}.
\]
\end{Lem}
\begin{pf}
Let's start with an alternate representation of the Legendre function
from \cite{T85}, page 158,
and a specific derivative
\begin{eqnarray*}
P_{-{1/2}+i t}(\cosh r) &=& \frac{\sqrt{2}}{\pi}
\int_0^r \frac{\cos(tu) \,du}{\sqrt{\cosh r - \cosh u}} ,\\
A(\cosh r) :\!&=&- \frac{\partial^2}{(\partial t)^2}
{ \bigg|_{t=0}}P_{-{1/2}+i t}(\cosh r)\\
&=& \frac{\sqrt{2}}{\pi}
\int_0^r\frac{u^2 \,du}{\sqrt{\cosh r - \cosh u}}.
\end{eqnarray*}
In the next step we make use of \cite{T85}, Exercise 28(b), page 158.
For fixed $t,\kappa>0 $, and large $\delta$, that is, $r=\kappa
/\delta
\to0$
\begin{eqnarray*}
P_{-{1/2}+it} \biggl(\cosh\frac{\kappa}{\delta} \biggr)
&\asymp& J_0 \biggl(t \frac{\kappa}{\delta} \biggr)\\
&=& \sum_{j=0}^\infty\frac{(-1)^j}{(j!)^2} \biggl(\frac{t\kappa
}{2\delta} \biggr)^{2j}\\
&=& 1 - \frac{1}{4} \biggl(\frac{t}{\delta} \biggr)^2 \kappa^2 + \cdots.
\end{eqnarray*}
Similarly,
\[
P_{-{1/2}+i {t}/{\delta}}(\cosh\kappa)
= P_{-{1/2}}(\cosh\kappa) - \frac{A(\kappa)}{2} \biggl(\frac
{t}{\delta} \biggr)^2 + \cdots.
\]
Then with the two above developments,
\begin{eqnarray*}
\mathcal{H} {H^{\delta}} \biggl({\frac1 2}+it \biggr)
&=& \frac{2\pi}{\delta}\int_0^\infty H(e^{-\kappa}i)
P_{-{1/2}}(\cosh\kappa) \biggl(1+O \biggl(\frac{1}{\delta^2} \biggr)
\biggr) \sinh(\kappa) \,d\kappa
\\
&\asymp&\frac{2\pi}{\delta}\int_0^\infty H(e^{-\kappa}i)
P_{-{1/2}+i {t}/{\delta}}(\cosh\kappa) \sinh(\kappa)\,
d\kappa
\\
&=& \frac{1}{\delta} \mathcal{H} H \biggl({\frac1 2}+i \frac{t}{\delta
} \biggr).
\end{eqnarray*}

Moreover since $\mathcal{H} H ({\frac1 2}+iu ) =0$ for $u \in(0,
{\frac1 2})$,
\begin{eqnarray*}
\|H^{\delta}\|^2_2 &\asymp&\frac{1}{4\pi}
\int_0^{\infty} \frac{1}{\delta^2}
\biggl|\mathcal{H} H \biggl({\frac1 2}+i \frac{t}{\delta} \biggr) \biggr|^2 t
\tanh(\pi t) \,dt\\
&=& \frac{1}{4\pi} \int_0^{\infty}
\biggl|\mathcal{H} H \biggl({\frac1 2}+iu \biggr) \biggr|^2
\frac{\tanh(\pi\delta u)}{\tanh(\pi u)} u \tanh(\pi u) \,du\\
&\asymp& \|H\|^2_2
\end{eqnarray*}
and taking additionally into account that $\delta$ is large
\begin{eqnarray*}
&& \| \Delta^{\alpha/2} H^\delta\|_2^2
\\
&&\qquad\asymp
\frac{1}{4\pi} \int_0^{\infty}
(t^2 + 1/4)^{\alpha}
\frac{1}{\delta^2}
\biggl|\mathcal{H} H \biggl({\frac1 2}+i \frac{t}{\delta} \biggr) \biggr|^2 t
\tanh(\pi t) \,dt\\
&&\qquad\le
\delta^{2\alpha} \frac{1}{4\pi} \int_0^{\infty}
\biggl( u^2 + \frac1 {4} \biggr)^{\alpha}
\biggl|\mathcal{H} H \biggl({\frac1 2}+iu \biggr) \biggr|^2
\frac{\tanh(\pi\delta u)}{\tanh(\pi u)} u \tanh(\pi u)\, du\\
&&\qquad\le C \delta^{2\alpha} \| \Delta^{\alpha/2} H \|_2^2 .
\end{eqnarray*}

Since by hypothesis
\[
|\mathcal{H} f_\varepsilon(1/2 + it) |
\le C_2 \exp[-|1/2 + it|^\beta]
= C_2 e^{-(1/4 + t^2)^{\beta/2}},
\]
we obtain
%
%
\begin{equation}\label{eq:Boundfhat}
\sup_{u \in[{1/ 2}, 1]} |\mathcal{H} f_\varepsilon(1/2 + i \delta u) |
\le C_2 e^{-(1/4 + \delta^2/4)^{\beta/2}}
\le C_2 e^{- ({\delta}/{2} )^{\beta}} .
\end{equation}
It follows from (\ref{eq:BoundedSupportofH}), (\ref{eq:DilateH(t)}) and
(\ref{eq:Boundfhat})
that there is a constant $C>0$ such that
\begin{eqnarray*}
\int| f_\varepsilon* H^\delta|^2
&=& \int| \mathcal{H} {H^\delta} (1/2 + it)|^2
|\mathcal{H} f_\varepsilon(1/2 + it)|^2 \,d\tau\\
&\asymp& \int_{1/2}^1 \delta^{-2}
|\mathcal{H} H (1/2 + iu)|^2 |\mathcal{H} f_\varepsilon(1/2 + i
\delta u)|^2
(\delta u) \tanh(\pi\delta u) \delta \,du \\
&\le& C e^{-2 ({\delta}/{2} )^{\beta}} .
\end{eqnarray*}
This completes the proof of Lemma \ref{lem:Dilate}.
\end{pf}

\subsubsection{Bound on $g_0 = f_\varepsilon*f_0$}\label{f0:scn}

Choose $f_0$ as
\[
f_0 (k(e^{-r}i) ) = f_0(e^{-r}i) = \frac{a-1} {2\pi} (\cosh r)^{-a}
\qquad \mbox{for }a > 1.
\]
\begin{Lem}
$f_0$ is $\mathbb{S}\mathbb{O}(2)$-invariant and
\[
\int_{r=0}^\infty\int_{u=0}^{2\pi} f_0 (k_u(e^{-r}i) ) \sinh r
\,dr \,du = 1.
\]
\end{Lem}
\begin{pf}
By the definition of $f_0$, it is $\mathbb{S}\mathbb{O}(2)$-invariant and
\begin{eqnarray*}
&&\int_{r=0}^\infty\int_{u=0}^{2\pi} f_0 (k_u(e^{-r}i) )
\sinh r \,dr \,du \\
&&\qquad= (a-1)
\int_{r=0}^\infty(\cosh r)^{-a} \sinh r \,dr \\
&&\qquad= (a-1) \int_{x=1}^\infty x^{-a} \,dx =1 .
\end{eqnarray*}
\upqed\end{pf}
\begin{Lem}\label{lem:Ceps}
Suppose that $f_\varepsilon$ is bi-invariant and
\[
\int_0^{C_\varepsilon} f_{\varepsilon}(e^{-s}i) \sinh s \,ds\geq\frac{1}{2}
\]
for some positive constant $C_\varepsilon$.
Then, $g_0 = f_\varepsilon* f_0$ satisfies
\[
g_0 (e^{-r} i) \ge\cases{
2\pi(a-1) 2^{a-1} e^{-2ar}, &\quad
for $r > C_\varepsilon$,\cr
2\pi(a-1) 2^{a-1} e^{-2aC_{\varepsilon}},&\quad
for $r \leq C_\varepsilon$.}
\]
\end{Lem}
\begin{pf}
Note that $0 \le s \le C_\varepsilon< r$ implies $r + s \le2r$
whereas $0 \le s,r \le C_\varepsilon$, implies $r + s \le
2C_\varepsilon$.
It follows from these facts, $f_{\varepsilon}\geq0$, Lemmas \ref
{lem:conv.rep} and \ref{lem:R} that
\begin{eqnarray*}
&&
f_{\varepsilon}*f_0(e^{-r}i) \\
&&\qquad= 2\pi\int_0^{2\pi}\int_0^\infty f_{\varepsilon}(e^{-s}i)
f_0 \bigl(k_{\eta(r,s,\phi)} e^{-R(r,s,\phi)}i \bigr) \sinh s \,ds\,
d\phi\\
&&\qquad= (a-1)\int_0^{2\pi}\int_0^\infty f_{\varepsilon}(e^{-s}i)
\cosh^{-a} (R(r,s,\phi) ) \sinh s \,ds \,d\phi\\
&&\qquad\geq2\pi(a-1) \int_0^{C_{\varepsilon}} f_{\varepsilon}(e^{-s}i)
\cosh^{-a} (r+s) \sinh s \,ds \\
&&\qquad\ge\frac{a-1}{\cosh^{a} (2\tau)} \int_0^{C_\varepsilon}
f_{\varepsilon}(e^{-s}i)
\sinh s \,ds \\
&&\qquad\ge2\pi\frac{a-1}{2 \cosh^{a} (2\tau)} \ge2\pi(a-1) 2^{a-1}
e^{-2a\tau}
\end{eqnarray*}
with $\tau= r$ for $r>C_{\varepsilon}$ and $\tau= C_{\varepsilon}$
for $r\leq C_{\varepsilon}$.
\end{pf}

\subsubsection{Bound on $f_\varepsilon* H^\delta$}\label{feps:scn}

We note from \cite{L72}, page 188,
\[
P_{-{1}/{2}} (\cosh r) = \frac{2}{\pi\cosh(r/2)}K \bigl( \tanh
(r/2) \bigr),
\]
where the complete elliptic integral of the first kind is defined by
\[
K(t) = \int_{0}^{\pi/2} \frac{d\phi} {\sqrt{1-t^2 \sin^2 \phi}} .
\]
In consequence there is a $C>0$ such that for all $r>0$
%
%
\begin{equation}\label{ineq:leb}
P_{- {1}/{2}}(\cosh r) \leq C .
\end{equation}

Define
\[
\mu_\delta(e^{-r}i) = \delta e^{-(m_0 - 1) \delta r}
\qquad \mbox{for }r \ge0.
\]
When $m_0 > 1$ and $\delta> 0$,
$\| \mu_\delta\|_\infty\le\delta$
and $\mu_\delta(e^{-r}i)$ is monotonically decreasing in~$r$.
\begin{Lem}\label{lem:m0}
Suppose $ | H(e^{-r}i) | \le C e^{-m_0 r}$ with $m_0 > 1$
and $\delta\ge1$.
Then, there is a constant $C>0$ such that $H^\delta(e^{-r}i) \le C \mu
_\delta(e^{-r}i)$.
\end{Lem}
\begin{pf}
By Taylor's expansion,
\[
(e^r)^\delta- (e^{-r})^\delta= \delta\eta^{\delta-1} (e^r - e^{-r})
\]
for $ e^{-r} \leq\eta\leq e^r $.
Since $\eta^{\delta-1} \leq e^{(\delta-1) r} \le e^{\delta r}$,
\[
\frac{\sinh(\delta r)}{\sinh r} \le\delta e^{\delta r}\qquad \mbox{for
}r \ge0.
\]
Hence with $C>0$ from (\ref{ineq:leb}),
\begin{eqnarray*}
| H^\delta(e^{-r}i) |
&=& | H(e^{-\delta r}i) |
\frac{\sinh(\delta r)}{\sinh r} P_{- 1 /2} (\cosh(\delta r)
)
\\
&\le& C | H(e^{-\delta r}i) | \frac{\sinh(\delta r)}{\sinh r}
\le C e^{-m_0 \delta r} \delta e^{\delta r}
= C \mu_\delta(e^{-r}i).
\end{eqnarray*}
\upqed\end{pf}
\begin{Lem}\label{lem:H}
Suppose that
\[
\int_{r - \xi_0 r}^{r + \xi_0 r} f_\varepsilon(e^{-s} i) \sinh s \,ds
\le C e^{-(\xi- \xi_0)r}
\qquad \mbox{for }0 < \xi_0 < 1 \mbox{ and } \xi> 1 + \xi_0,
\]
$H(e^{-r} i)$ is bounded and monotonely decreasing in $r$
and satisfies the tail condition
%
%
\begin{equation}\label{cond:TailH}
H(e^{-r}i) \le C e^{-m_0 r}
\qquad \mbox{for }m_0 \xi_0 > \xi.
\end{equation}
Then, for $\delta\ge1$ and $r > 0$,
\[
(f_\varepsilon* H^\delta)(e^{-r}i) \le C \delta e^{-(\xi- \xi_0)r}.
\]
\end{Lem}
\begin{pf}
Set $\mathcal{R}_1 = \{s \dvtx|r-s| \le\xi_0 r \}$ and
$\mathcal{R}_2 = \{s \dvtx |r-s| > \xi_0 r \}$.
It follows from Lemmas \ref{lem:conv.rep}, \ref{lem:b.m.H} and
\ref{lem:m0} that
%
%
\begin{eqnarray}\label{eq:f*Hdelta}
&&\frac{1}{2\pi} (f_\varepsilon* H^\delta)(e^{-r}i)\nonumber\\
&&\qquad= \int_{\phi=0}^{2\pi} \int_{s=0}^\infty
f_\varepsilon(e^{-s} i) H^\delta\bigl(k_{\eta(r,s,\phi)} e^{-R(r,s,\phi
)} i \bigr) \sinh s \,d\phi \,ds
\nonumber\\
&&\qquad\le C \int_{\phi=0}^{2\pi} \int_{s=0}^\infty
f_\varepsilon(e^{-s} i) \mu_\delta\bigl(e^{-R(r,s,\phi)} i \bigr) \sinh
s \,d\phi \,ds
\nonumber\\
&&\qquad\le
2\pi C \int f_\varepsilon(e^{-s} i) \mu_\delta\bigl(e^{-|r-s|} i \bigr)
\sinh s \,ds.
\end{eqnarray}
Set
\[
I_j = \int_{\mathcal{R}_j} f_\varepsilon(e^{-s} i) \mu_\delta
\bigl(e^{-|r-s|} i \bigr) \sinh s \,ds
\qquad \mbox{for }j = 1, 2.
\]
Since $m_0 > \xi/\xi_0 > 1$,
\[
\| \mu_\delta\|_\infty\le\delta.
\]
Observe that
%
%
\begin{eqnarray}\label{eq:I1}
I_1
&\le& \| \mu_\delta\|_\infty\int_{|r-s| \le\xi_0 r}
f_\varepsilon
(e^{-s} i) \sinh s \,ds
\nonumber\\
&\le& C\delta e^{-(\xi- \xi_0)r}.
\end{eqnarray}
From the tail condition (\ref{cond:TailH}),
%
%
\begin{eqnarray}\label{eq:I2}
I_2
&=& \int_{|r-s| > \xi_0 r} f_\varepsilon(e^{-s} i) \bigl( \delta
e^{-(m_0 - 1) \delta|r-s|} \bigr)
\sinh s \,ds
\nonumber\\
&\le& \delta e^{-(m_0 - 1) \delta\xi_0 r} \int_{|r-s| > \xi_0 r}
f_\varepsilon(e^{-s} i) \sinh s \,ds
\nonumber\\
&\le& \delta e^{-(m_0 - 1) \delta\xi_0 r}.
\end{eqnarray}
Since $r \ge0$, $m_0 \xi_0 > \xi$ and $\delta$ is large,
\begin{eqnarray*}
\frac{e^{-(\xi- \xi_0)r}} {e^{-(m_0 - 1) \delta\xi_0 r}}
&=& e^{(m_0 - 1) \delta\xi_0 r - (\xi- \xi_0)r}
\ge e^{(m_0 -1 ) \xi_0 r - (\xi- \xi_0)r}
= e^{(m_0 \xi_0 - \xi)r}
\ge1.
\end{eqnarray*}
Combining (\ref{eq:f*Hdelta}), (\ref{eq:I1}), (\ref{eq:I2}), we have
the desired result.
\end{pf}

\subsubsection{\texorpdfstring{Chi-square distance: Proof of
(\protect\ref{eq:chisq-1})}{Chi-square distance: Proof of (B.6)}}\label{chi2:dist:scn}

With the above notation, choose a pair of densities
\[
f_0 (k(e^{-r}i) ) = f_0(e^{-r}i) = \frac{a-1} {2\pi} (\cosh
r)^{-a} \quad \mbox{and} \quad
f_n = f_0 + C_H \delta^{-\alpha} H^\delta,
\]
where $a$ satisfies $1 < a < 2$, and $H$ satisfies the hypotheses of
Lemmas \ref{lem:Dilate} and~\ref{lem:H}.
By choosing $C_H$ close to 0, we have $f_0, f_n \in{\mathcal
{F}_{\alpha}(Q)}$ for all large
$\delta$.
Let $g_0 = f_\varepsilon* f_0$ and $g_n = f_\varepsilon* f_n$ with
$f_{\varepsilon}$ satisfying the hypotheses of Lemma \ref{lem:H} and
(D.4) with $\gamma=1$.

The $\chi^2$ distance between $g_0 $ and $g_n$ is defined by
\begin{eqnarray*}
\chi^2(g_0, g_n) := \int\frac{ (g_n - g_0 )^2} {g_0} \,dz =
2\pi C_H^2 \delta^{-2\alpha}
\int\frac{ (f_\varepsilon* H^\delta(e^{-r}i) )^2}
{g_0(e^{-r}i)} \sinh r \,dr
\end{eqnarray*}
and with a suitable constant $M>0$ guaranteed by Lemma \ref{lem:Ceps},
by Lemmas \ref{lem:Dilate}, \ref{lem:m0} and \ref{lem:H}
\begin{eqnarray*}
&&\int\frac{ (f_\varepsilon* H^\delta(e^{-r}i) )^2}
{g_0(e^{-r}i)} \sinh r \,dr
\\
&&\qquad\le\biggl(\int_{e^r \le M} + \int_{e^r > M} \biggr)
\frac{ (f_\varepsilon* H^\delta(e^{-r}i) )^2} {g_0(e^{-r}i)}
\sinh r \,dr
\\
&&\qquad\le
\frac{M^{2a}} {C} \int\bigl(f_\varepsilon* H^\delta(e^{-r}i) \bigr)^2
\sinh r \,dr
+ \frac{\delta^2}{C} \int_{e^r > M}
\frac{e^{-2(\xi- \xi_0)r}} {e^{-2ar}} \sinh r \,dr
\\
&&\qquad= O \bigl( M^{2a} e^{-2 ({\delta}/{2} )^{\beta}} + \delta
^2 M^{-2(\xi- \xi_0) + 2a + 1} \bigr)
= O \bigl( e^{-\mu_1 ({\delta}/{2} )^{\beta}} \bigr),
\end{eqnarray*}
where $C = 2\pi(a-1)2^{a-1}$. For the last equality we set
$M = e^{{1}/{2} ({\delta}/{2} )^{\beta}}$,
$\mu_0 = 2(\xi- \xi_0) - 2a - 1$
and
$\mu_1 = \min(\mu_0/2, (2-a) )$.
Then indeed
\[
M^{2a} e^{-2 ({\delta}/{2} )^{\beta}} = e^{-(2-a)
({\delta}/{2} )^{\beta}} = O ( e^{-\mu_1 \delta^{\beta
}} )
\quad \mbox{and} \quad
\delta^2 M^{-\mu_0}
= O ( e^{-\mu_1 \delta^{\beta}} ).
\]
Hence,
%
%
\begin{equation}\label{eq:chi2Bound}
\chi^2(g_0, g_n) = O \bigl( \delta^{ -2\alpha} e^{-\mu_1 (
{\delta}/{2} )^{\beta}} \bigr).
\end{equation}
Letting $e^{-\mu_1 ({\delta}/{2} )^{\beta}} = n^{-1}$, or
equivalently
%
%
\begin{equation}\label{eq:delta-rate}
\delta= 2\mu_1^{-1/ {\beta}} ( \log n )^{1/ {\beta}},
\end{equation}
we conclude that the right-hand side of (\ref{eq:chi2Bound}) is of
order $o(n^{-1})$, that is, (\ref{eq:chisq-1}) is proven.
\end{appendix}

\section*{Acknowledgments}
The authors would like to express their sincerest gratitude to the late
Dennis M. Healy Jr. (1957--2009), whose untimely death was most
unfortunate. The second author had been in conversations with him
about this very problem for over twenty years, and he has passed away
prior to seeing the completion of this remarkable problem. For this
and so much more, we are all grateful.

\printaddresses


\begin{thebibliography}{99}

\bibitem{AS98}
\textsc{Abramovich, F.} and \textsc{Silverman, B.} (1998).
{Wavelet decomposition approaches to statistical inverse problems}.
\textit{Biometrika} \textbf{85} 115--129.
\MR{1627226}

\bibitem{AH03}
\textsc{Allen, J. C.} and \textsc{Healy Jr., D. M.} (2003).
{Hyperbolic geometry, {N}ehari's theorem, electric circuits, and
analog signal processing}.
\textit{Modern Signal Processing} \textbf{46} 1--62.

\bibitem{BP03}
\textsc{Bhattacharya, R.} and \textsc{Patrangenaru, V.} (2003).
Large sample theory of intrinsic and extrinsic sample means on
manifolds {I}.
\textit{Ann. Statist.} \textbf{31} 1--29.
\MR{1962498}

\bibitem{BP05}
\textsc{Bhattacharya, R.} and \textsc{Patrangenaru, V.} (2005).
Large sample theory of intrinsic and extrinsic sample means on
manifolds {II}.
\textit{Ann. Statist.} \textbf{33} 1225--1259.
\MR{2195634}

\bibitem{BHMR07}
\textsc{Bissantz, N.}, \textsc{Hohage, T.}, \textsc{Munk, A.} and
\textsc{Ruymgaart, F.} (2007).
Convergence rates of general regularization methods for statistical
inverse problems and applications.
\textit{SIAM J. Numer. Anal.} \textbf{45} 2610--2636.
\MR{2361904}

\bibitem{B01}
\textsc{Brown, B. H.} (2001).
Medical impedance tomography and process impedance tomography: A
brief review.
\textit{Meas. Sci. Technol.} \textbf{12} 991--996.

\bibitem{BT07}
\textsc{Butucea, C.} and \textsc{Tsybakov, A.} (2007).
Sharp optimality for density deconvolution with dominating bias, {I};
{II}.
\textit{Theory Probab. Appl.} \textbf{52} 111--128; 336--349.
\MR{2354572}

\bibitem{CGLT03}
\textsc{Cavalier, L.}, \textsc{Golubev, G.}, \textsc{Lepski, O.} and
\textsc{Tsybakov, A.} (2003).
Block thresholding and sharp adaptive estimation in severely
ill-posed inverse problems.
\textit{Theory Probab. Appl.} \textbf{48} 534--556.
\MR{2141349}

\bibitem{CGPT02}
\textsc{Cavalier, L.}, \textsc{Golubev, G.}, \textsc{Picard, D.} and
\textsc{Tsybakov, A.} (2002).
Oracle inequalities for inverse problems.
\textit{Ann. Statist.} \textbf{30} 843--874.
\MR{1922543}

\bibitem{Collin92}
\textsc{Collin, R. E.} (1992).
\textit{Foundations for {M}icrowave {E}ngineering}.
McGraw-Hill, New York.

\bibitem{DG04}
\textsc{Delaigle, A.} and \textsc{Gijbels, I.} (2004).
{Practical bandwidth selection in deconvolution kernel density
estimation.}
\textit{Comput. Statist. Data Anal.} \textbf{45} 249--267.
\MR{2045631}

\bibitem{DHM08}
\textsc{Delaigle, A.}, \textsc{Hall, P.} and \textsc{Meister, A.}
(2008).
{On deconvolution with repeated measurements.}
\textit{Ann. Statist.} \textbf{36} 665--685.
\MR{2396811}

\bibitem{DMR96}
\textsc{Dey, A.}, \textsc{Mair, B.} and \textsc{Ruymgaart, F.}
(1996).
Cross-validation for parameter selection in inverse estimation
problems.
\textit{Scand. J. Statist.} \textbf{23} 609--620.
\MR{1439715}

\bibitem{EMCG02}
\textsc{Ewart, G.}, \textsc{Mills, K.}, \textsc{Cox, G.} and
\textsc{Gage, P.} (2002).
Amiloride derivatives block ion channel activity and enhancement of
virus-like particle budding caused by {HIV}-1 protein {V}pu.
\textit{Eur. Biophys. J.} \textbf{31} 26--35.

\bibitem{Fan91}
\textsc{Fan, J.} (1991).
On the optimal rates of convergence for nonparametric deconvolution
problem.
\textit{Ann. Statist.} \textbf{19} 1257--1272.
\MR{1126324}

\bibitem{Gold99}
\textsc{Goldenshluger, A.} (1999).
On pointwise adaptive nonparametric deconvolution.
\textit{Bernoulli} \textbf{5} 907--926.
\MR{1715444}

\bibitem{HJDH08}
\textsc{Hahn, G.}, \textsc{Just, A.}, \textsc{Dittmar, J.} and
\textsc{Hellig, G.} (2008).
Systematic errors of {EIT} systems determined by easily-scalable
resistive phantoms.
\textit{Physiol. Meas.} \textbf{29} 163--172.

\bibitem{HM07}
\textsc{Hall, P.} and \textsc{Meister, A.} (2007).
A ridge-parameter approach to deconvolution.
\textit{Ann. Statist.} \textbf{35} 1535--1558.
\MR{2351096}

\bibitem{HQ05}
\textsc{Hall, P.} and \textsc{Qiu, P.} (2005).
{Discrete-transform approach to deconvolution problems.}
\textit{Biometrika} \textbf{92} 135--148.
\MR{2158615}

\bibitem{H82}
\textsc{Helton, J. W.} (1982).
Non-{E}uclidean functional analysis and electronics.
\textit{Bull. Amer. Math. Soc. (N.S.)} \textbf{7} 1--64.
\MR{0656197}

\bibitem{JKPR04}
\textsc{Johnstone, I. M.}, \textsc{Kerkyacharian, G.}, \textsc
{Picard, D.} and \textsc{Raimondo, M.} (2004).
{Wavelet deconvolution in a periodic setting.}
\textit{J. R. Stat. Soc. Ser. B Stat. Methodol.} \textbf{66} 547--573.
\MR{2088290}

\bibitem{JR04}
\textsc{Johnstone, I. M.} and \textsc{Raimondo, M.} (2004).
{Periodic boxcar deconvolution and Diophantine approximation.}
\textit{Ann. Statist.} \textbf{32} 1781--1804.
\MR{2102493}

\bibitem{JoSi91}
\textsc{Johnstone, I. M.} and \textsc{Silverman, B. W.} (1991).
Discretization effects in statistical inverse problems.
\textit{J. Complexity} \textbf{7} 1--34.
\MR{1096170}

\bibitem{KM03}
\textsc{Kalifa, J.} and \textsc{Mallat, S.} (2003).
{Thresholding estimators for linear inverse problems and
deconvolutions.}
\textit{Ann. Statist.} \textbf{31} 58--109.
\MR{1962500}

\bibitem{KV97}
\textsc{Kass, R. E.} and \textsc{Vos, P. W.} (1997).
\textit{Geometrical {F}oundations of {A}symptotic {I}nference}.
Wiley, New York.
\MR{1461540}

\bibitem{KPP07}
\textsc{Kerkyacharian, G.}, \textsc{Petrushev, P.}, \textsc{Picard,
D. B.} and \textsc{Willer, T.} (2007).
Needlet algorithms for estimation in inverse problems.
\textit{Electron. J. Stat.} \textbf{1} 30--76.
\MR{2312145}

\bibitem{KK05}
\textsc{Kim, P. T.} and \textsc{Koo, J.-Y.} (2005).
Statistical inverse problems on manifolds.
\textit{J. Fourier Anal. Appl.} \textbf{11}
639--653.
\MR{2190676}

\bibitem{KR01}
\textsc{Kim, P. T.} and \textsc{Richards, D. S.} (2001).
Deconvolution density estimation on compact {L}ie groups.
In \textit{Algebraic Methods in Statistics and Probability ({N}otre
{D}ame, {IN}, 2000)}. \textit{Contemporary Mathematics}
\textbf{287} 155--171. Amer. Math. Soc.,
Providence, RI.
\MR{1873674}


\bibitem{KR02}
\textsc{Kim, P. T.} and \textsc{Richards, D. S.} (2008).
Diffusion tensor imaging and deconvolution on the space of positive
definite symmetric matrices.
In \textit{Mathematical Foundations of Computational Anatomy: Geometric
and Statistical Methods for Biological Shape Variability Modeling}
(X. Pennec and S. Joshi, eds.) 140--149.
Available at
\href{http://www-sop.inria.fr/asclepios/events/MFCA08/Proceedings/MFCA08_Proceedings.pdf}{http://www-sop.inria.fr/asclepios/events/MFCA08/Proceedings/}
\href{http://www-sop.inria.fr/asclepios/events/MFCA08/Proceedings/MFCA08_Proceedings.pdf}{MFCA08\_Proceedings.pdf}.

\bibitem{L72}
\textsc{Lebedev, N.} (1972).
\textit{Special Functions and Their Applications}.
Dover, New York.
\MR{0350075}

\bibitem{Lee97}
\textsc{Lee, J. M.} (1997).
\textit{Riemannian {M}anifolds: {A}n {I}ntroduction to {C}urvature}.
Springer, New York.
\MR{1468735}

\bibitem{LMR00}
\textsc{Lovric, M.}, \textsc{Min-Oo, M.} and \textsc{Ruh, E. A.}
(2000).
Multivariate normal distributions parametrized as a riemannian
symmetric space.
\textit{J. Multivariate Anal.} \textbf{74} 36--48.
\MR{1790612}

\bibitem{MR96}
\textsc{Mair, B. A.} and \textsc{Ruymgaart, F. H.} (1996).
Statistical inverse estimation in {H}ilbert scales.
\textit{SIAM J. Appl. Math.} \textbf{56} 1424--1444.
\MR{1409127}

\bibitem{M92}
\textsc{McCullaugh, P.} (1992).
Conditional inference and cauchy models.
\textit{Biometrika} \textbf{79} 247--259.
\MR{1185127}

\bibitem{M96}
\textsc{McCullaugh, P.} (1996).
M\"{o}bius transformation and {C}auchy parameter estimation.
\textit{Ann. Statist.} \textbf{24} 787--808.
\MR{1394988}

\bibitem{MM04}
\textsc{Mizera, I.} and \textsc{M\"uller, C.} (2004).
Location-scale depth (with discussion).
\textit{J. Amer. Statist. Assoc.} \textbf{99} 949--966.
\MR{2109488}

\bibitem{N08}
\textsc{Neubauer, A.} (2008).
The convergence of a new heuristic parameter selection criterion for
general regularization methods.
\textit{Inverse Problems} \textbf{24} 055005.
\MR{2438940}

\bibitem{NePa64}
\textsc{Nevanlinna, R.} and \textsc{Paatero, V.} (1964).
\textit{Einf\"uhrung in die Funtkionentheorie}.
Birkh\"auser, Basel.

\bibitem{PS09}
\textsc{Pensky, M.} and \textsc{Sapatinas, T.} (2009).
{Functional deconvolution in a periodic setting: Uniform case.}
\textit{Ann. Statist.} \textbf{37} 73--104.
\MR{2488345}

\bibitem{PS05}
\textsc{Pereverzev, S.} and \textsc{Schock, E.} (2005).
{On the adaptive selection of the parameter in regularization of
ill-posed problems.}
\textit{SIAM J. Numer. Anal.} \textbf{43} 2060--2076.
\MR{2192331}

\bibitem{RS04}
\textsc{R\"omer, W.} and \textsc{Steinem, C.} (2004).
Impedance analysis and single-channel recordings on nano-black lipid
membranes based on porous alumina.
\textit{Biophys. J.} \textbf{86} 955--965.

\bibitem{S03}
\textsc{Szkutnik, Z.} (2003).
{Doubly smoothed EM algorithm for statistical inverse problems.}
\textit{J.~Amer. Statist. Assoc.} \textbf{98} 178--190.
\MR{1965684}

\bibitem{T85}
\textsc{Terras, A.} (1985).
\textit{Harmonic Analysis on Symmetric Spaces and Applications.} I.
Springer, New York.
\MR{0791406}

\bibitem{VEGS08}
\textsc{Van Es, B.}, \textsc{Gugushvili, S.} and \textsc
{Spreij, P.}
(2008).
{Deconvolution for an atomic distribution}.
\textit{Electron. J. Stat.} \textbf{2} 265--297.
\MR{2399196}

\bibitem{W94}
\textsc{Wand, M. P.} and \textsc{Jones, M. C.} (1994).
\textit{Kernel Smoothing}.
Chapman and Hall, London.
\MR{1319818}

\end{thebibliography}
\end{document}